\documentclass[11pt]{article}
\usepackage{graphicx}
\usepackage{amsmath,amsfonts,amssymb}
\usepackage{epstopdf}
\usepackage{url}
\usepackage{amsmath}
\usepackage{graphicx}
\usepackage{epsfig}
\usepackage{color}
\def\tto{\;{\lower 1pt \hbox{$\rightarrow$}}\kern -10pt
\hbox{\raise 2pt \hbox{$\rightarrow$}}\;}
\def\Hat{\widehat}

\def\ra{\rangle}
\def\la{\langle}

\def\B{I\!\!B}

\def\ox{\bar{x}}

\def\cone{\mbox{\rm cone}\,}

\def\O{\Omega}
\def\ph{\varphi}

\newcounter{lk}

\begin{document}

\begin{center}
\vspace*{0.3in} \textbf{VARIATIONAL ANALYSIS OF DIRECTIONAL MINIMAL TIME
FUNCTIONS AND APPLICATIONS TO LOCATION PROBLEMS}\\[2ex]Nguyen Mau
Nam\footnote{Fariborz Maseeh Department of Mathematics and Statistics,
Portland State University, Portland, OR 97202, United States (email:
mau.nam.nguyen@pdx.edu). The research of Nguyen Mau Nam was partially
supported by the Simons Foundation under grant \#208785.} and Constantin
Z\u{a}linescu\footnote{University Al.I. Cuza Ia\c{s}i, Faculty of Mathematics,
700506 Ia\c{s}i, Romania (email: zalinesc@uaic.ro).}\\[2ex]
\end{center}

{\small \textbf{Abstract.} This paper is devoted to the study of
\emph{directional minimal time functions} that specify the minimal time for a
vector to reach an object following its given direction. We provide a careful
analysis of general and generalized differentiation properties of this class
of functions. The analysis allows us to study a new model of facility location
that involves sets. This is a continuation of our effort in applying
variational analysis to facility location problems. }

\medskip
\vspace*{0,05in} \noindent {\bf Key words.} Directional minimal time functions, scalarization functions, generalized \\differentiation, facility location problems.

{\small \newtheorem{Theorem}{Theorem}[section]
\newtheorem{Proposition}[Theorem]{Proposition}
\newtheorem{Remark}[Theorem]{Remark} \newtheorem{Lemma}[Theorem]{Lemma}
\newtheorem{Corollary}[Theorem]{Corollary}
\newtheorem{Definition}[Theorem]{Definition}
\newtheorem{Example}[Theorem]{Example}
\renewcommand{\theequation}{\thesection.\arabic{equation}} }

\section{Introduction and Preliminaries}

Let $X$ be a real normed linear space. Given a vector $v\in X$, $v\neq0$,
and a nonempty closed set $\Omega\subseteq X$, the directional minimal time
function with direction $v$ and target set $\Omega$ is defined by
\begin{equation}
T_{v}(x;\Omega):=\inf\{t\geq0\mid x+tv\in\Omega\}. \label{s}%
\end{equation}
This class of functions is similar to the class of
scalarization functions that has been used to study vector optimization problems:
\begin{equation}\varphi_{v}(x;\Omega):=\inf\{t\in\mathbb{R}\mid
x+tv\in\Omega\};\label{phi}\end{equation} see \cite{ti,tz} and the
references therein. We will see later on that if $v\in\Omega_{\infty}$, then \begin{equation*}
T_{v}(x;\Omega)=\max\left\{  \varphi_{v}(x;\Omega),0\right\}
\quad\forall x\in X. \label{t-phi}\end{equation*}
Lipschitz properties and a formula for computing
subdifferentials in the sense of convex analysis of scalarization functions
were studied in \cite{tz}, but generalized differentiation properties
involving nonconvex structures have not been considered in the literature.
Notice that the directional minimal time function (\ref{s}) is a particular
case of general minimal time functions considered, e.g., in \cite{cowo,bmn11} and the references therein. However,
the specific structure of the function makes it distinct from the general
case.\vspace*{0.05in}

In this paper, we mainly study generalized differentiation properties the directional minimal time function
(\ref{s}) and the scalarization function (\ref{phi}), as well as applications to facility location problems. The location model
that motivates our study is a generalized version of the celebrated
Fermat-Torricelli problem: given a finite number of nonempty closed target
sets $\Omega_{i}$ for $i=1,\ldots,n$ and $n$ nonzero vectors $v_{i}$ for
$i=1,\ldots,n$, along with a nonempty closed constraint set $\Omega_0$, find a
point $\bar{x}\in\Omega_0$ to place the initial points of the vectors such that
the total times for the vectors to reach the target sets is minimal. This
problem can be modeled as the following optimization problem:
\begin{equation}
\mbox{\rm minimize }\sum_{i=1}^{n}T_{v_{i}}(x;\Omega_{i}%
)\mbox{ subject to }x\in\Omega_0. \label{ft}%
\end{equation}
The location model of this type seems to be very interesting, but it has not
been considered in the literature even in the convex case. Since the functions
involved in the problem are not differentiable in general, our approach
involves developing new tools of modern variational analysis for these
functions and apply them to solve the problem. Let us recall basic definitions
and properties of variational analysis that will play an important role in the
sequel.\vspace*{0.05in}

Let $\Omega\subseteq X$ and let $\bar{x}\in\Omega$. A vector $x^{\ast}\in
X^{\ast}$ is called a Fr\'{e}chet normal to $\Omega$ at $\bar{x}$ if
\[
\langle x^{\ast},x-\bar{x}\rangle\leq\mbox{\rm o}(\Vert x-\bar{x}%
\Vert)\mbox{ for }x\in\Omega.
\]
The set of all Fr\'{e}chet normals to $\Omega$ at $\bar{x}$ is called the
Fr\'{e}chet normal cone to $\Omega$ at $\bar{x}$, denoted by $\widehat{N}%
(\bar{x};\Omega)$. \vspace*{0.05in}

A vector $x^{\ast}\in X^{\ast}$ is called a Mordukhovich/limiting normal to
$\Omega$ at $\bar{x}$ if there are sequences $x_{k}\xrightarrow{\O}\bar{x}$
and $x_{k}^{\ast}\xrightarrow{w^*}x^{\ast}$ with $x_{k}^{\ast}\in\widehat
{N}(x_{k};\Omega)$ for every $k$. The set of all Mordukhovich normals to
$\Omega$ at $\bar{x}$ is called the Mordukhovich/limiting normal cone to the
set at this point. In this definition, the notation $x_{k}\xrightarrow{\O}\bar
{x}$ means that $x_{k}\rightarrow\bar{x}$ and $x_{k}\in\Omega$ for every $k$.
\vspace*{0.05in}

Let $\psi:X\rightarrow(-\infty,\infty]$ be an extended real-valued function
and let $\bar{x}$ be an element of the domain of the function
$\mbox{\rm dom }\psi:=\{x\in X\mid\psi(x)<\infty\}$. The Fr\'{e}chet
subdifferential of $\psi$ at $\bar{x}$ is defined by
\[
\widehat{\partial}\psi(\bar{x}):=\{x^{\ast}\in X^{\ast}\mid\langle x^{\ast
},x-\bar{x}\rangle\leq\psi(x)-\psi(\bar{x})+\mbox{\rm o}(\Vert x-\bar{x}%
\Vert)\} \mbox{ for }x\in X.
\]
The Mordukhovich/limiting subdifferential of $\psi$ at $\bar{x}$, denoted by
$\partial\psi(\bar{x})$, is the set of all vectors $x^{\ast}\in X^{\ast}$ such
that there exist sequences $x_{k}\xrightarrow{\psi}\bar{x}$, and $x_{k}^{\ast
}\in\widehat{\partial}\psi(x_{k})$ with $x_{k}^{\ast}\xrightarrow{w^*}x^{\ast
}$, where $x_{k}\xrightarrow{\psi}\bar{x}$ means that $x_{k}\rightarrow\bar
{x}$ and $\psi(x_{k})\rightarrow\psi(\bar{x})$.

In the case where the set $\Omega$ is convex or the function $\psi$ is convex,
the corresponding Fr\'{e}chet normal cone and subdifferential structures reduce respectively to those in the sense of convex analysis. If $X$ is a Banach space, then the same holds for the limiting structures.

The singular subdifferential of $\psi$ at $\bar{x}$, denoted by $\partial
^{\infty}\psi(\bar{x})$, is the set of all vectors $x^{\ast}\in X^{\ast}$ such
that there exist sequences $\lambda_{k}\downarrow0$, $x_{k}%
\xrightarrow{\psi}\bar{x}$, and $x_{k}^{\ast}\in\widehat{\partial}\psi(x_{k})$
with $\lambda_{k}x_{k}^{\ast}\xrightarrow{w^*}x^{\ast}$, where $x_{k}%
\xrightarrow{\psi}\bar{x}$ means that $x_{k}\rightarrow\bar{x}$ and
$\psi(x_{k})\rightarrow\psi(\bar{x})$, and $\lambda_{k}\downarrow0$ means that
$\lambda_{k}\rightarrow0$ and $\lambda_{k}>0$ for every $k$. This
subdifferential structure is particularly important in recognizing the
Lipschitz continuity of functions; see, e.g., \cite[Theorem 3.52]{mor} and
\cite[Theorem 9.13]{rw}.\vspace*{0.05in}

The paper is organized as follows. Section 2 presents general properties of
the directional minimal time function (\ref{s}). In Section 3, we study
generalized differentiation properties of the function that involve various
kinds of subdifferential structures in convex and nonconvex settings. Although
some results in Section 2 and 3 can be derived from \cite{bmn11}, we provide
detailed simplified proofs for the convenience of the reader. Section 4 is
devoted to the study of Lipschitzian properties using both direct and
generalized differentiation approaches. Finally, in Section 5, we apply the
results from the previous sections to study location problems (\ref{ft}%
).\vspace*{0.05in}

Throughout the paper, we use the following standing assumptions unless
otherwise stated: $\Omega$ is a nonempty closed subset of a real normed linear
space $X$; $v$ is a nonzero vector in $X$. Moreover, when there is
no risk of confusion, we will use $T$ instead of
$T_{v}(\cdot,\Omega)$ and $\varphi$ instead of
$\varphi_{v}(\cdot,\Omega).$ For a nonzero vector $v\in X$, we will use the following notations:
\begin{align*}
&\cone\{v\}:=\{\lambda v\mid \lambda\geq 0\}=\Bbb R_{+}v,\\
&\mbox{\rm span }\{v\}:=\{\lambda v\mid \lambda\in \Bbb R\}=\Bbb R v,\\
&\{v\}^+:=\{x^*\in X^*\mid \la x^*, v\ra \geq 0\},\\
&\{v\}^\perp:=\{x^*\in X^*\mid \la x^*, v\ra =0\}.
\end{align*}
\section{General Properties}

Let us start with simple representations of the domain and the epigraph of the directional minimal time function (\ref{s}). Recall that the recession
cone of $\Omega$ is given by
\[
\Omega_{\infty}:=\{u\in X\mid\omega+\lambda u\in\Omega
\;\mbox{\rm for all }\omega\in\Omega\;\mbox{\rm and for all }\lambda\geq0\}.
\]

\begin{Proposition}\label{domain} The domain of the directional minimal
time function {\rm(\ref{s})}  is given by
\begin{equation*}
\mbox{\rm dom }T=\O-\cone\{v\}.
\end{equation*}
Suppose further that $v\in \O_\infty$. Then
\begin{gather}\label{d2}
\mbox{\rm dom }T=\O-\mbox{\rm span }\{v\},\\
\label{epi} \mbox{\rm epi }T=\{(x,t)\in X\times{\Bbb
R}\mid  t\geq 0, x+tv\in \O\},
\end{gather}
and \begin{equation} T(x)=\max\left\{
\varphi(x),0\right\} \quad\forall x\in X.
\label{t-phi}\end{equation}
\end{Proposition}
Proof: By the definition,
\begin{align*}
\mbox{\rm dom }T &  =\{x\in X\mid T(x)<\infty\}\\
&  =\{x\in X\mid x+tv\in\Omega\mbox{ for some }t\geq0\}\\
&  =\{x\in X\mid x\in\Omega-tv\mbox{ for some }t\geq0\}\\
&  =\Omega-\mbox{\rm cone}\,\{v\}.
\end{align*}

Assume that $v\in\Omega_{\infty}$. For any $x\in\Omega-\mbox{\rm span }\{v\}$,
one has $x=\omega-\lambda v$, where $\omega\in\Omega$ and $\lambda
\in\mathbb{R}$. If $\lambda\geq0$, then $x\in\Omega
-\mbox{\rm cone }\{v\}=\mbox{\rm dom }T$. In the case
$\lambda<0$, since $v\in\Omega_{\infty}$, one has $x=\omega+(-\lambda
)v\in\Omega\subseteq\mbox{\rm dom
}T$. Thus, (\ref{d2}) holds in this case. The proof for
(\ref{epi}) is straightforward.

The inequality $\geq$ in (\ref{t-phi}) (even for
arbitrary $v$) is obvious. Let $\max\left\{  \varphi(x),0\right\}
<\lambda$. By the definition of $\varphi$, there exists
$t\in(-\infty,\lambda)$ such that $x+tv\in\Omega$. Then $x+\lambda
v=x+tv+(\lambda-t)v\in\Omega+\mathbb{R}_{+}v=\Omega$, and so
$T(x)\leq\lambda$. It follows that $T(x)\leq\max\left\{
\varphi(x),0\right\} .$ $\hfill\square$

\begin{Remark}{\rm The \emph{scalarization
function} associated with $\O$ and $v$ was introduced by Gerstewitz
(Tammer) and Iwanow \cite{ti} as in (\ref{phi}).
This function has been used extensively in vector optimization; see
\cite{tz} and the references therein. Two important properties of
$\varphi_{v}(\cdot;\Omega)$ are:
\begin{equation}
\varphi_v(x+tv;\O)=\varphi_v(x;\O)-t \mbox{ for all }x\in X,\ t\in
{\Bbb R},\label{p-phi}
\end{equation}
and
\begin{equation*}
\varphi_{v}(\cdot;\Omega)=\varphi_{v}(\cdot;\Omega+{\Bbb R}_{+}v).
\end{equation*}
However, $T_v(\cdot;\O)\geq T_v(\cdot; \O+{\Bbb R}_{+}v)$ in
general. This can be seen by taking $X={\Bbb R}^2$,
$\O=[-1,1]\times[-1,1]$, and $v=(1,0)$.}
\end{Remark}

\begin{Proposition}\label{bdry} The infimum in the directional minimal
time function {\rm (\ref{s})} always attains for any $x\in \mbox{\rm
dom } T$. That means
\begin{equation*}
\Pi_v(x;\O):=x+T(x)v\in \O
\end{equation*}
for all $x\in \mbox{\rm dom }T$. In fact, $\Pi_v(x;\O)
\in\mbox{\rm bd }\O$ for all $x$ with $T(x)\in (0,\infty)$.
\end{Proposition}Proof: Let $x\in\mbox{dom }T$ and let
$t:=T(x)$. Then there exists a sequence $t_{k}\rightarrow t$ with
$t_{k}\in\lbrack0,\infty)$ and $x+t_{k}v\in\Omega$ for every $k$. Thus,
$x+tv\in\Omega$ since $\Omega$ is closed. Suppose $t\in(0,\infty)$ and
$\Pi_v(x;\Omega)\notin\mbox{\rm bd }\Omega$. Then $\Pi_v(x;\Omega)\in
\mbox{\rm int }\Omega$, so there exists $\delta>0$ such that
$x+tv+I\!\!B(0;\delta)\subseteq\Omega$ and $\dfrac{\delta}{\Vert v\Vert}<t$.
Then $x+tv-\dfrac{\delta}{\Vert v\Vert}v\in\Omega$. This implies
$t=T(x)\leq t-\dfrac{\delta}{\Vert v\Vert}$, which is a
contradiction. $\hfill\square$ \vspace*{0.05in}

For any $x\in\mbox{\rm dom }T$, the element $\Pi
_{v}(x;\Omega)=x+T(x)v$ is called the projection from $x$ to
$\Omega$ with respect to the directional minimal time function (\ref{s}).

\begin{Proposition}\label{lsc} For any $\alpha\in\Bbb R$, define
$$\mathcal{L}_\alpha:=\{x\in X\mid  T(x)\leq \alpha\} \;\mbox{\rm and }
\mathcal{L}^<_\alpha:=\{x\in X\mid  T(x)< \alpha\}.$$
Then
\begin{equation}\label{in1}
\mathcal{L}_\alpha=\O-[0, \alpha]v \mbox{ and }
\mathcal{L}^{<}_\alpha=\O-[0, \alpha)v\; \mbox{\rm for }\alpha\geq 0.
\end{equation}
The equality $T(x)=0$ holds if and only if $x\in \O$. Moreover,
$T$ is lower semicontinuous.
\end{Proposition}Proof: The proof of (\ref{in1}) is obvious. Suppose
$x\in\Omega$. Then $x+0v\in\Omega$, so $T(x)\leq0$, and hence
$T(x)=0$. Now suppose $T(x)=0$. Then $x+T(x)v=x\in\Omega$. Since $\mathcal{L}_{\alpha}$ is closed for every $\alpha\in \Bbb R$
($\mathcal{L}_{\alpha}$ being empty for $\alpha<0$), the function $T$
is lower semicontinuous. $\hfill\square$

\begin{Remark}\label{rem2.5}
{\rm Because $T$ is lsc, in the case
$v\in \O _\infty$, by (\ref{t-phi}), we have that $T$
and $\ph$ coincide on a neighborhood of $x$ for any
$x\in X$ with $T(x)>0$. This shows that the properties of
$T$ proved at elements $x\in X$ with $T(x)>0$ can
be extended to similar properties of  $\ph$ [having in
view also (\ref{p-phi})].}
\end{Remark}

\begin{Proposition} The function $T$ is convex if and only if
$\O$ is convex.
\end{Proposition}Proof: Suppose that $\Omega$ is convex. Fix $x_{1},x_{2}%
\in\mbox{\rm dom }T$ and $\lambda\in(0,1)$. Then
\[
x_{1}+T(x_{1})v\in\Omega\mbox{ and }x_{2}+T(x_{2})v\in\Omega.
\]
Since $\Omega$ is convex, this implies
\[
\lambda x_{1}+(1-\lambda)x_{2}+[\lambda T(x_{1})+(1-\lambda
)T(x_{2})]v\in\Omega.
\]
Thus,
\[
T(\lambda x_{1}+(1-\lambda)x_{2})\leq\lambda T(x_{1})+(1-\lambda)T(x_{2}).
\]
Conversely, if $T$ is convex, then $\{x\in X\mid
T(x)\leq0\}=\Omega$ is convex. $\hfill\square$ \vspace*{0.05in}

We are now going to study the strict convexity of $T$. For
$a,b\in X$, let
\[
\lbrack a,b]:=\{ta+(1-t)b\mid t\in\lbrack0,1]\}
\]
be the line segment connecting the two points.

A set $\Omega$ is called strictly convex if for any $x_{1},x_{2}\in\Omega$,
$x_{1}\neq x_{2}$, and for any $t\in(0,1)$, one has $tx_{1}+(1-t)x_{2}%
\in\mbox{\rm int }\Omega$.

\begin{Proposition}\label{2.5} Suppose that $[a,b]\subseteq \mbox{\rm dom }
T\setminus \O$ and $b-a\notin \mbox{\rm span}\{v\}$. If
$\O$ is strictly convex, then $T$ is strictly convex on
$[a,b]$.
\end{Proposition}Proof: Suppose by contradiction that $x_{1},x_{2}\in\lbrack
a,b]$, $x_{1}\neq x_{2}$, $\alpha\in(0,1)$ and
\[
T(\alpha x_{1}+(1-\alpha)x_{2})=\alpha T(x_{1})+(1-\alpha)T(x_{2})>0.
\]
Then
\begin{align*}
\Pi_{v}(\alpha x_{1}+(1-\alpha)x_{2};\Omega)  &  =\alpha x_{1}+(1-\alpha
)x_{2}+T(\alpha x_{1}+(1-\alpha)x_{2})v\\
&  =\alpha(x_{1}+T(x_{1})v)+(1-\alpha)(x_{2}+T(x_{2})v)\\
&  =\alpha\Pi_{v}(x_{1};\Omega)+(1-\alpha)\Pi_{v}(x_{2};\Omega).
\end{align*}
Since $x_{1},x_{2}\in\mbox{\rm dom }T$ and $x_{1}%
-x_{2}\notin \Bbb Rv$, one has that $\Pi_{v}(x_{1};\Omega)\neq\Pi_{v}(x_{2}%
;\Omega)$. Consequently,
\[
\alpha\Pi_{v}(x_{1};\Omega)+(1-\alpha)\Pi_{v}(x_{2};\Omega)\in\mbox{\rm int
}\Omega,
\]
which is a contradiction since the projection must belong to the boundary of
$\Omega$ by \newline Proposition \ref{bdry}. $\hfill\square$

\begin{Proposition}\label{time property} Fix $\ox\in\mbox{\rm dom }
T$. One always has
\begin{equation}\label{t1}
T(\ox+\lambda v) = T(\ox)-\lambda \mbox{ for any
}0\leq\lambda\leq T(\ox).
\end{equation}
For $\lambda>0$, suppose further that $\ox-\gamma v\notin\O$ for
every $\gamma\in (0, \lambda]$. Then
$$T(\ox-\lambda v)=T(\ox)+\lambda.$$
Moreover, \begin{equation}\label{maxin}
T(x+tv)\geq \max\{T(x)-t, 0\} \mbox{ for all }x\in X
\mbox{ and }t\geq 0,
\end{equation}
and the equality holds if $v\in\O_\infty$.
\end{Proposition}
Proof: Since $\bar{x}+T(\bar{x})v\in\Omega$, one has $\bar
{x}+\lambda v+(T(\bar{x})-\lambda)v\in\Omega$. Thus, $$T(\bar
{x}+\lambda v)\leq T(\bar{x})-\lambda<\infty.$$ Let us show
that
\begin{equation}
T(\bar{x})\leq\lambda+T(\bar{x}+\lambda v).
\label{lm2.8}%
\end{equation}
Indeed, one can assume $t:=T(\ox+\lambda v)<\infty$. Then
\[
\ox+\lambda v+tv\in\Omega.
\]
This implies $T(\ox)\leq\lambda+t$, and the result follows. From
(\ref{lm2.8}), one has $T(\bar{x})-\lambda\leq T(\bar
{x}+\lambda v)$, and hence (\ref{t1}) holds.

Let us now prove the second equality under the assumption that $\bar{x}-\gamma
v\notin\Omega$ for every $\gamma\in(0,\lambda]$. Let $\bar{y}:=\bar{x}-\lambda
v$. Then $\bar{x}=\bar{y}+\lambda v$. From (\ref{lm2.8}), one has
\[
T(\bar{y})=T(\bar{x}-\lambda v)\leq T(\bar{x})+\lambda<\infty.
\]
This implies $\bar{y}\in\mbox{\rm dom }T$. In this case, we
can easily see that $0<\lambda\leq T(\bar{y})$. Indeed, let
$\eta:=T(\bar{y})$. Then $\bar{y}+\eta v=\bar{x}-\lambda v+\eta
v=\bar{x}-(\lambda-\eta)v\in\Omega$. This implies $\lambda-\eta\notin
(0,\lambda]$, so $\lambda\leq\eta=T(\bar{y})$. Applying (\ref{t1}),
one has $T(\bar{x})=T(\bar{y})-\lambda=T(\bar
{x}-\lambda v)-\lambda$, and the second equality follows.

The inequality (\ref{maxin}) follows from (\ref{lm2.8}). Let us prove that the
equality holds if $v\in \Omega_{\infty}$. There is nothing to prove if
$T(x)=\infty$. Consider two cases $T(x)>t$ and
$T(x)\leq t$. In the first case, one has
\[
x+tv+(T(x)-t)v=x+T(x)v\in\Omega.
\]
Thus, $T(x+tv)\leq T(x)-t$, and the conclusion holds. In
the second case, one has $x+tv=x+T(x)v+(t-T(x))v\in\Omega+\O_\infty=\Omega$. Thus, $T(x+tv)=0$, and the conclusion also
holds. $\hfill\square$

\begin{Corollary}\label{cor2.9}Assume that there exist $\bar{x}\in\Omega$
and $\bar{t}>0$ such that $\bar{x}+\bar{t}v\in\Omega$
but $[\bar{x},\bar {x}+\bar{t}v]\not \subseteq \Omega$.
Then there exists $t\in\lbrack 0,\bar{t})$ such that
$T$ is finite and not continuous at
$\bar{x}+tv.$
\end{Corollary}
Proof: By hypothesis, there exists
$t_{0}\in(0,\bar{t})$ such that
$\bar{x}+t_{0}v\notin\Omega$. Considering the largest interval
$I\subseteq\lbrack0,\bar{t}]$ containing $t_{0}\in I$ such that
$\bar{x}+sv\notin\Omega$ for every $s\in I$, we may assume that
$I=(0,\bar{t})$. Thus, for every $\lambda\in(0,\bar{t})$, we
have that $\bar{x}+\bar{t}v-\gamma v\notin\Omega$ for
every $\gamma \in(0,\lambda]$. Using the previous proposition, we get
$T(\bar {x}+\bar{t}v-\lambda
v)=T(\bar{x}+\bar{t}v)+\lambda=\lambda$, and so
$$\lim_{t\rightarrow0_{+}}T(\bar{x}+tv)=\lim_{\lambda
\rightarrow\bar{t}_{-}}T(\bar{x}+\bar{t}v-\lambda
v)=\bar{t}\neq T(\bar{x})=0.$$
Therefore, $T$ is not
continuous at $\bar{x}$.  $\hfill\square$

The next result provides a sufficient condition and
a necessary condition for the continuity of $T$ at
some $\ox\in\mbox{\rm dom }T$.

\begin{Proposition}\label{p2.8} Let $\ox\in\mbox{\rm dom }T$.
If $\ox + (T(\ox), T(\ox)+\gamma)v\subseteq \mbox{\rm int
}\O$ for some $\gamma>0$, then $T$ is continuous at
$\ox$. Suppose additionally that  $v\in\O_\infty$. Then
$T$ is continuous at $\ox$ if and only if $\ox +
(T(\ox), \infty)v\subseteq \mbox{\rm int }\O$.
\end{Proposition}Proof: Fix any number $\lambda$ such that $T(\bar
{x})<\lambda$. We can assume that $\lambda<T(\bar{x})+\gamma$.
Then $\bar{x}+\lambda v\in\mbox{\rm int }\Omega$. Choose $\delta>0$ such that
\[
\bar{x}+\lambda v+\delta I\!\!B\subseteq\Omega.
\]
For any $x\in\bar{x}+\delta I\!\!B$, one has $T(x)\leq\lambda$.
Thus, $T$ is upper semicontinuous at $\bar{x}$, so it is
continuous at this point.

Suppose that $v\in\Omega_{\infty}$ and $T$ is continuous at
$\bar{x}$. Consider $\lambda>T(\bar{x})$. Since $T$ is continuous at $\ox$, there exists $\delta>0$ such that whenever $x\in
V:=\bar{x}+\delta I\!\!B$, one has $T(x)<\lambda$. This implies
$x+\lambda v\in\Omega$, so $x\in\Omega-\lambda v$. Thus, $V\subseteq
\Omega-\lambda v$ or $\lambda v+V=\bar{x}+\lambda v+\delta I\!\!B\subseteq
\Omega$. Therefore, $\bar{x}+\lambda v\in\mbox{\rm int
}\Omega$.$\hfill\square$

\begin{Example} Let $\O:=\{(x,y)\in \Bbb R^2\mid  y=
0\}$ and let $v=(0,1)$. Then $T$ is continuous at
$\ox=(0,-1)$, but $\ox + (T(\ox), \infty)v\nsubseteq \mbox{\rm
int }\O$.
\end{Example}

We end this section with other properties of the
directional minimal time function which will be used in the next
sections.

\begin{Proposition}\label{ll1} For two nonempty closed subsets $A$ and
$B$ of $X$, the following hold:\\[1ex]
{\rm (1)} If $A\subseteq B$, then $T_v(x, B)\leq T_v(x, A)$ for all $x\in X$.\\
{\rm (2)} $T_v(x+y; A+B)\leq T_v(x; A)+T_v(y; B)$ for all $x,y\in X$.
\end{Proposition}Proof: (1) It is obvious if $T_{v}(x;A)=\infty$. Suppose
$t:=T_{v}(x;A)<\infty$. Then
\[
x+tv\in A\subseteq B.
\]
Thus, $T_{v}(x;B)\leq t=T_{v}(x;A)$.\newline(2) The conclusion if obvious if
$T_{v}(x;A)=\infty$ or $T_{v}(y;B)=\infty$. In the other case, let
$s:=T_{v}(x;A)$ and $t:=T_{v}(x;B)$. Then
\[
x+tv\in A\mbox{ and }y+tv\in B.
\]
Thus, $x+y+(s+t)v\in A+B$. This implies
\[
T_{v}(x+y;A+B)\leq s+t=T_{v}(x;A)+T_{v}(y;B).
\]
The proof is now complete. $\hfill\square$

\section{Generalized Differentiation Properties}

In this section, we are going to study generalized differentiation properties
of the directional minimal time function (\ref{s}). Various subdifferential
structures of variational analysis will be employed to study the function. The
results from this section will be important for the study of Lipschitz
continuity of the function in Section 4 and for applications to location
problems in Section 5.

\subsection{Fenchel conjugate and sub\-gradients in the sense of convex
analysis}

For a function $\psi:X\rightarrow(-\infty,\infty]$, recall that the Fenchel
conjugate of $\psi$ is an extended real-valued function on $X^{\ast}$ defined
by
\[
\psi^{\ast}(x^{\ast}):=\sup\{\langle x^{\ast},x\rangle-\psi(x)\mid x\in X\}.
\]
Let us start with a formula for representing the Fenchel conjugate of the
directional minimal time function (\ref{s}) in terms of the support function of $\O$ defined on $X^*$ by
\[
\sigma_\O(x^*):=\sup\{\la x^*,\omega\ra \mid \omega\in \O\}.
\]
\begin{Proposition}\label{fenchel} The function $T$ is a lsc proper function and
\begin{equation*}
T^*(x^*)=\begin{cases}
\sigma_\O(x^*), &\text{if }\; \la x^*,-v\ra\leq 1, \\
\infty, & \text{otherwise}.
\end{cases}
\end{equation*}
\end{Proposition}Proof: The fact that $T$ is lower semicontinuous has been
proved in Proposition \ref{lsc}. It is proper since $\Omega\subseteq
\mbox{\rm dom }T$. We have the following
\begin{align*}
T^{\ast}(x^{\ast})  &  =\sup\{\langle x^{\ast},x\rangle-T(x)\mid
x\in\mbox{\rm dom }T\}\\
&  =\sup\{\langle x^{\ast},x\rangle-t\mid t\geq0,x+tv\in\Omega\}\\
&  =\sup\{\langle x^{\ast},\omega\rangle+t[\langle x^{\ast},-v\rangle-1]\mid
t\geq0,\omega\in\Omega\}\\
&  =\sup\{\langle x^{\ast},\omega\rangle\mid\omega\in\Omega\}+\sup\{t(\langle
x^{\ast},-v\rangle-1)\mid t\geq0\}.
\end{align*}
The formula then follows easily. $\hfill\square$

\begin{Theorem}\label{cv1} Let
$\ox\in \mbox{\rm dom }T$, where $\O$ is convex. Then
\begin{equation}
\partial T(\ox)=\left\{
\begin{array}
[c]{ll}%
N(\ox,\Omega)\cap\left\{  x^{\ast}\in
X^{\ast}\mid\left\langle
v,x^{\ast}\right\rangle \geq-1\right\}, \text{if }\ox\in\Omega,\\
N(\ox+T(\ox)v,\Omega)\cap\left\{  x^{\ast}\in X^{\ast}%
\mid\left\langle v,x^{\ast}\right\rangle =-1\right\}, \text{otherwise.}%
\end{array}
\right.  \label{r-subd}%
\end{equation}
Moreover, if $\ox\in\O$ and $v\in\O_\infty$, then \begin{equation}
\label{c2}\partial T(\ox)=\{x^*\in X^*\;
|\; -1\leq \la x^*, v\ra \leq 0\}\cap N(\ox;\O).\end{equation}
\end{Theorem}Proof: Consider $x^{\ast}\in\partial T(\bar{x})$. Then $T(\bar
{x})+T^{\ast}(x^{\ast})=\left\langle \bar{x},x^{\ast}\right\rangle
$. By Proposition \ref{fenchel},
\[
\left\langle v,x^{\ast}\right\rangle \geq-1\mbox{ and }T(\bar{x}%
)+\sigma_{\Omega}(x^{\ast})=\left\langle \bar{x},x^{\ast}\right\rangle .
\]
Since $\bar{x}+T(\bar{x})v\in\Omega$, the following holds
\[
\left\langle \bar{x},x^{\ast}\right\rangle \geq T(\bar{x})+\left\langle
\bar{x}+T(\bar{x})v,x^{\ast}\right\rangle .
\]
It follows that
\[
T(\bar{x})[1+\left\langle v,x^{\ast}\right\rangle ]\leq0.
\]
Because $T(\bar{x})\geq0$ and $1+\left\langle v,x^{\ast}\right\rangle \geq0,$
one has
\[
T(\bar{x})[1+\left\langle v,x^{\ast}\right\rangle ]=0.
\]
If $\bar{x}\in\Omega$, then $\sigma_{\Omega}(x^{\ast})=\left\langle
\bar {x},x^{\ast}\right\rangle $, and so $x^{\ast}\in
N(\bar{x},\Omega)$. If $\bar{x}\notin\Omega$, then
$T(\bar{x})>0$, and so $1+\left\langle v,x^{\ast }\right\rangle =0.$
It follows that
\[
\sigma_{\Omega}(x^{\ast})=\left\langle \bar{x}+T(\bar{x})v,x^{\ast
}\right\rangle ,
\]
and $x^{\ast}\in N(\bar{x}+T(\bar{x})v,\Omega)$. Thus, the
inclusion $\subseteq$ holds in (\ref{r-subd}). Conversely, if
$\bar{x}\in\Omega$ and $x^{\ast}\in N(\bar{x},\Omega)\cap\left\{
x^{\ast}\in X^{\ast}\mid \left\langle v,x^{\ast}\right\rangle
\geq-1\right\}  $, then
\[
T(\bar{x})[1+\left\langle v,x^{\ast}\right\rangle ]=0\mbox{ and
}\sigma_{\Omega}(x^{\ast})=\left\langle \bar{x},x^{\ast}\right\rangle .
\]
Thus, $T(\bar{x})+T^{\ast}(x^{\ast})=\left\langle \bar{x},x^{\ast
}\right\rangle $, and so $x^{\ast}\in\partial T(\bar{x})$. If $\bar{x}%
\notin\Omega$ and $x^{\ast}\in N(\bar{x}+T(\bar{x})v,\Omega)\cap\left\{
x^{\ast}\in X^{\ast}\mid\left\langle v,x^{\ast}\right\rangle =-1\right\}  $,
then
\[
T(\bar{x})[1+\left\langle v,x^{\ast}\right\rangle ]=0\mbox{
and }\sigma_{\Omega}(x^{\ast})=\left\langle \bar{x}+T(\bar{x})v,x^{\ast
}\right\rangle .
\]
It follows that $T(\bar{x})+T^{\ast}(x^{\ast})=\left\langle
\bar {x},x^{\ast}\right\rangle $, and again $x^{\ast}\in\partial
T(\bar{x}).$

Under the condition $v\in\Omega_{\infty}$, one sees easily that $\langle
x^{\ast},v\rangle\leq0$ for every $x^{\ast}\in N(\bar{x};\Omega)$. Thus,
equality (\ref{c2}) follows. $\hfill\square$

\begin{Proposition}\label{cv2} Let $\ox\in \O$, where $\O$ is convex and $X$ is a Banach space. Then
\begin{equation*}
\partial^\infty T(\ox)=N(\ox,\Omega)\cap \{v\}^+=N(\ox,\operatorname*{dom}T).
\end{equation*}
Moreover,
if $v\in\Omega_{\infty}$, then \begin{equation*}
\partial^\infty T(\ox)=N(\ox,\Omega)\cap \{v\}^\perp=N(\ox,
\operatorname*{dom}T).
\end{equation*}
\end{Proposition}Proof: Fix
any $x^{*}\in\partial^{\infty}T(\bar{x})$. Then there exist sequences $x_{k}%
\xrightarrow{T}\bar{x}$, $\lambda_{k}\downarrow0$, $\lambda_{k}x_{k}%
^{*}\xrightarrow{w^*}x^{*}$ such that $x_{k}^{*}\in\partial T(x_{k})$. Let
$\tilde{x}_{k}:=x_{k}+T(x_{k})v\in\Omega$. By Theorem \ref{cv1}, $x_{k}^{*}\in
N(\tilde{x}_{k};\Omega)$ and $\langle x_{k}^{*},-v\rangle\leq1$ (the equality
holds if $x_{k}\notin\Omega$). Then $\lambda_kx^*_k\in N(\tilde{x}_{k};\Omega)$, and $(\lambda_kx^*_k)$ is a bounded sequence as $X$ is a Banach space. Since $\tilde{x}_{k}:=x_{k}+T(x_{k}%
)v\rightarrow\bar{x}$, one has $x^{*}\in N(\bar{x};\Omega)$. Moreover,
$\langle\lambda_{k}x_{k}^{*},-v\rangle\leq\lambda_{k}$, and hence $\langle
x^{*},-v\rangle\leq0$. It follows that $x^{*}\in\{v\}^{+}$. Now fix any
$x^{*}\in N(\bar{x},\Omega)\cap\{v\}^{+}$. Then $kx^{*}\in N(\bar{x};\Omega)$
and $\langle kx^{*},-v\rangle\leq0<1$ for every $k$. By Theorem \ref{cv1},
$kx^{*}\in\partial T(\bar{x})$, and hence $x^{*}\in\dfrac{1}{k}\partial
T(\bar{x})$. By definition, $x^{*}\in\partial^{\infty}T(\bar{x})$.

We have seen that $\operatorname*{dom}T=\Omega-\Bbb R_{+}v$, and so
\begin{align*}
x^{\ast}\in N(\bar{x},\operatorname*{dom}T)  &  \Leftrightarrow\left\langle
x-tv-\bar{x},x^{\ast}\right\rangle \leq0\ \forall x\in\Omega,\ t\geq0\\
&  \Leftrightarrow\left\langle x-\bar{x},x^{\ast}\right\rangle \leq0\ \forall
x\in\Omega\text{ and }\left\langle v,x^{\ast}\right\rangle \geq0.
\end{align*}
This implies $N(\bar{x},\operatorname*{dom}T)=N(\bar{x},\Omega)\cap\{v\}^{+}$.

The proof for the second equalities under the condition $v\in\Omega_{\infty}$
follows from the last observation in the proof of Theorem \ref{cv1}.
$\hfill\square$

\begin{Proposition} Let $\bar{x}%
\in\operatorname*{dom}T\setminus\Omega$, where $\O$ is convex and $X$ is a Banach space. Then
$$\partial^{\infty}T(\bar {x})\subseteq
N(\bar{x}+T(\bar{x})v,\Omega)\cap\{v\}^{\perp}.$$
Moreover, if $\partial T(\bar{x})$ is nonempty, then equality
holds.
\end{Proposition}Proof: Let $x^{\ast}\in\partial^{\infty}T(\bar{x}).$
Then there exist sequences $x_{k}\overset{T}{\rightarrow}\bar{x},$
$\lambda _{k}\downarrow0$,
$\lambda_{k}x_{k}^{\ast}\overset{w^{\ast}}{\rightarrow }x^{*}$
such that $x_{k}^{\ast}\in\partial T(x_{k})$ for $k\geq1$. Since
$T(x_{k})\rightarrow T(\bar{x})>0$, we may assume that
$T(x_{k})>0$ for $k\geq1$. Let
$\tilde{x}_{k}:=x_{k}+T(x_{k})v\in\Omega$. Clearly,
$\tilde{x}_{k}\rightarrow\tilde{x}:=\bar{x}+T(\bar{x}%
)v\in\Omega$. By Theorem \ref{cv1},
\[
x_{k}^{\ast}\in N(\tilde{x}_{k};\Omega)\mbox{ and }\left\langle v,x_{k}^{\ast
}\right\rangle =-1.
\]
Because $X$ is a Banach space, we have that $(\lambda_{k}x_{k}^{\ast})$ is
bounded. From the inequality $\left\langle x-\tilde{x}_{k},x_{k}^{\ast
}\right\rangle \leq0$ for $k\geq1$ and $x\in\Omega$, one has
\[
\left\langle x-\tilde{x}_{k},\lambda_{k}x_{k}^{\ast}\right\rangle \leq0.
\]
This implies $\left\langle x-\tilde{x},x^{\ast}\right\rangle
\leq0$ for $x\in\Omega$. As $\left\langle
v,\lambda_{k}x_{k}^{\ast}\right\rangle =-\lambda_{k}\rightarrow0,$
we also have that $\left\langle v,x^{\ast }\right\rangle =0.$
Therefore, $x^{\ast}\in N(\tilde{x},\Omega )\cap\{v\}^{\perp}$.

Assume now that $\partial T(\bar{x})\neq\emptyset$ and fix $x_{0}^{\ast
}\in\partial T(\bar{x})$. By Theorem \ref{cv1}, we have that $x_{0}%
^{\ast}\in N(\bar{x}+T(\bar{x})v,\Omega)$ and $\left\langle
v,x_{0}^{\ast}\right\rangle =-1$. Take $x^{\ast}\in N(\bar{x}%
+T(\bar{x})v,\Omega)\cap\{v\}^{\perp}$. Then
\[
x_{k}^{\ast}:=x_{0}^{\ast}+kx^{\ast}\in N(\bar{x}+T(\bar{x}%
)v,\Omega)\mbox{ and }\left\langle v,x_{k}^{\ast}\right\rangle =-1.
\]
Thus, $x_{k}^{\ast}\in\partial T(\bar{x})$. Taking
$x_{k}:=\bar{x}$ and $\lambda_{k}:=1/k$, we obtain that
$x^{\ast}$ $\in\partial^{\infty }T(\bar{x})$. $\hfill\square$

\begin{Example}{\rm Let $\O=\B(0; \sqrt{8})$ in $\Bbb R^2$ with the
Euclidean norm and let $v=(1,1)$. For $\ox=(-2,-2)\in\O$, using the
formula from Theorem \ref{cv1}, one has
\begin{equation*}
\partial T(\ox)=\{(t, t)\mid  -1/2\leq t\leq 0\}.
\end{equation*}
For $\ox=(2, 2)$, one has
\begin{equation*}
\partial T(\ox)=\{(t,t)\mid  t\geq 0\}.
\end{equation*}
For $\ox=(-3,-3)\notin\O$, one has
\begin{equation*}
\partial T(\ox)=\{(-1/2, -1/2)\}.
\end{equation*}}
\end{Example}

We end this section with a result referring to the
scalarization function $\varphi$ given in (\ref{phi}). The result is
immediate from Proposition 3.4 and Theorem 3.2, using Remark
\ref{rem2.5}.

\begin{Corollary} Let $\bar{x}\in X$ satisfy
$\varphi(\bar{x})\in\mathbb{R}$, where $\O$ is convex. Assume that $v\in\Omega_{\infty}$. Then
\[
\partial\varphi(\bar{x})=N(\bar{x}+\varphi%
(\bar{x})v; \O)\cap\{x^{\ast}\in X^{\ast}\mid\left\langle x^{\ast
},v\right\rangle =-1\}.
\]
Moreover, if $X$ is a Banach space, then
\[
\partial^{\infty}\varphi(\bar{x})\subseteq N(\bar
{x}+\varphi(\bar{x})v; \O)\cap\{v\}^{\perp},
\]
and the equality holds if
$\partial\varphi(\bar{x})\neq\emptyset.$
\end{Corollary}

\subsection{Dini-Hadamard directional derivatives and sub\-gradients}

Let $\psi:X\rightarrow(-\infty,\infty]$ be a function and let $\bar{x}%
\in\mbox{\rm dom }\psi$. The Dini-Hadamard directional derivative of the
function at $\bar{x}$ in the direction $u$ is given by
\[
\psi^{\prime}(\bar{x},u):=\liminf_{t\rightarrow0^{+},\ v\rightarrow u}%
\frac{\psi(\bar{x}+tv)-\psi(\bar{x})}{t}.
\]
If $\psi$ is Lipschitz continuous around $\bar{x}$, the Dini-Hadamard
directional derivative coincides with the (lower) Dini directional derivative
\[
\psi^{\prime}(\bar{x};h):=\liminf_{t\rightarrow0^{+}}\dfrac{\psi(\bar
{x}+th)-\psi(\bar{x})}{t}.
\]
Let $\bar{x}\in\Omega$. The Bouligand contingent cone to $\Omega$ at $\bar{x}%
$, denoted by $K(\bar{x};\Omega)$, is the set of $d\in X$ such that there
exist sequences $t_{k}\rightarrow0^{+}$ and $d_{k}\rightarrow d$ with $\bar{x}%
+t_{k}d_{k}\in\Omega$ for every $k$. It is well-known that $K(\bar{x};\Omega)$
is a closed cone. Moreover, $K((x,\psi(x));\mbox{\rm epi}\,\psi
)=\mbox{\rm epi}\,\psi^{\prime}(x;\cdot)$ if $x\in\mbox{\rm dom }\psi$%
.\vspace*{0.05in}

We also define the Dini-Hadamard normal cone
\[
N^{-}(\bar{x};\Omega):=\{x^{\ast}\in X^{\ast}\mid\langle x^{\ast},d\rangle
\leq0\mbox{ for all }d\in K(\bar{x};\Omega)\},
\]
and the Dini-Hadamard subdifferential of $\psi$ at $\bar{x}$
\[
\partial^{-}\psi(\bar{x})=\{x^{\ast}\in X^{\ast}\mid\langle x^{\ast}%
,h\rangle\leq\psi^{\prime}(\bar{x};h)\mbox{ for all }h\in X\}.
\]
\begin{Lemma}\label{c1} Let $\ox\in\O$. Then
\begin{equation*}
T^\prime(\ox; h)=T_v(h; K(\ox; \O)).
\end{equation*}
\end{Lemma}Proof: Let $\bar{x}\in\Omega$. Set $\lambda:=T^{\prime
}(\bar{x},u);$ because $T(x)\geq T(\bar{x})=0$ for every
$x\in X$, we have that $\lambda\geq0$. If $\lambda=+\infty$, we have
clearly that $T_{v}(u,K(\bar{x},\Omega))\leq
T^{\prime}(\bar{x},u)$. Assume that $\lambda<\infty$. Then
there exist sequences $u_{k}\rightarrow u$ and
$t_{k}\rightarrow0^{+}$ such that $$\lambda_{k}:=t_{k}^{-1}\left[
T(\bar{x}+t_{k}u_{k})-T(\bar{x})\right]
=t_{k}^{-1}T(\bar {x}+t_{k}u_{k})\rightarrow\lambda.$$ In
particular, we may assume that
$\lambda_{k}\in \Bbb R_{+}$ for $k\geq1$. It follows that $\bar{x}+t_{k}%
(u_{k}+\lambda_{k}v)\in\Omega$, whence $u+\lambda v\in K(\bar{x}%
,\Omega)$, and so $T_{v}(u,K(\bar{x},\Omega))\leq\lambda$.
Therefore, $T_{v}(u,K(\bar{x},\Omega))\leq
T^{\prime}(\bar{x},u).$

Conversely, set
$\lambda:=T_{v}(u,K(\bar{x},\Omega))\in\lbrack0,\infty].$
Assume that $\lambda<\infty$. Then $u+\lambda v\in
K(\bar{x},\Omega)$, and so there exist sequences
$u_{k}^{\prime}\rightarrow u+\lambda v$ and
$t_{k}\rightarrow0^{+}$ such that $\bar{x}+t_{k}u_{k}^{\prime
}=\bar{x}+t_{k}(u_{k}^{\prime}-\lambda v)+t_{k}\lambda
v\in\Omega$ for $k\geq1$. Then
$T(\bar{x}+t_{k}(u_{k}^{\prime}-\lambda v))\leq t_{k}\lambda.$
Since $T(\bar{x})=0$ and $(u_{k}^{\prime}-\lambda v)\rightarrow
u$, we obtain that $T^{\prime}(\bar{x},u)\leq\lambda$. Thus,
$T^{\prime}(\bar{x},u)\leq T_{v}(u,K(\bar{x},\Omega))$ for
every $u\in X$. $\hfill\square$

\begin{Theorem}\label{DH1} For any $\ox\in\O$, one has
\begin{equation}\label{clr}
\partial^-T(\ox)=\{x^*\mid  \la x^*, -v\ra \leq 1\}\cap N^-(\ox;\O).
\end{equation}
\end{Theorem}Proof: Fix any $x^{*}\in\partial^{-}T(\bar{x})$. By the
definition and Lemma \ref{c1},
\[
\langle x^{*},h\rangle\leq T^{\prime}(\bar{x};h)=T_{v}(h;K(\bar{x}%
;\Omega))\mbox{ for all }h\in X.
\]
For $h=-v$, one has $\langle x^{*},-v\rangle\leq T_{v}(-v;K(\bar{x}%
;\Omega))\leq1$ since $0\in K(\bar{x};\Omega)$. Moreover, for all $h\in
K(\bar{x};\Omega)$, one has
\[
\langle x^{*},h\rangle\leq T_{v}(h;K(\bar{x};\Omega))=0.
\]
Thus, $x^{*}\in N^{-}(\bar{x};\Omega)$.

Now suppose that $\langle x^{*},-v\rangle\leq1$ and $x^{*}\in N^{-}(\bar
{x};\Omega)$. For any $h\in X$, let us show that
\[
\langle x^{*},h\rangle\leq T^{\prime}(\bar{x};h).
\]
The inequality holds obviously if $T^{\prime}(\bar{x};h)=\infty$. Consider the
case where
\[
t:=T_{v}(h;K(\bar{x};\Omega))<\infty.
\]
Then $h+tv\in K(\bar{x};\Omega)$. Thus, $\langle x^{*},h+tv\rangle\leq0$. This
implies
\[
\langle x^{*},h\rangle\leq t\langle x^{*},-v\rangle\leq t=T_{v}(h;K(\bar
{x};\Omega))=T^{\prime}(\bar{x};h).
\]
Therefore, $x^{*}\in\partial^{-}T(\bar{x})$, and (\ref{clr}) has been proved.
$\hfill\square$

\begin{Lemma}\label{lm3.8} Let $\bar{x}\in\operatorname*{dom}T
\setminus\Omega$ and let $\tilde{x}:=\Pi_v(\ox;\O)$. Then
\begin{equation*}
T^{\prime}(\bar{x},u)\leq\varphi_v (u,K(\tilde {x},\Omega))
\mbox{ for every }u\in X.
\end{equation*}
Moreover, equality holds for those $u\in X$ with
$T^{\prime}(\bar{x},u)\neq-\infty$. This equality also holds if
$v\in \O_\infty$.
\end{Lemma}
Proof: Assume that $\varphi_{v} (u,K(\tilde{x},\Omega)) <\lambda
\in\mathbb{R.}$ Then there exists $\mu\in(-\infty,\lambda)$ such
that $u+\mu v\in K(\tilde{x},\Omega)$, and so there exist
$(u_{k}^{\prime})\rightarrow
u+\mu v$ and $(t_{k})\rightarrow0^{+}$ such that $\tilde{x}+t_{k}%
u_{k}^{\prime}=\bar{x}+t_{k}(u_{k}^{\prime}-\mu v)+(T(\bar{x}%
)+t_{k}\mu)v\in\Omega$ for $k\geq1$. Since $T(\bar{x})>0,$
there exists $k_{0}\geq1$ such that $T(\bar{x})+t_{k}\mu>0$ for
$k\geq k_{0}$. It follows that
\[
T(\bar{x}+t_{k}(u_{k}^{\prime}-\mu v))\leq T(\bar{x})+t_{k}\mu.
\]
Since $u_{k}^{\prime}-\mu v\rightarrow u$, one has
$T^{\prime}(\bar {x},u)\leq\mu<\lambda$. Thus,
$T^{\prime}(\bar{x},u)\leq\varphi_v (u,K(\tilde {x},\Omega)).$

Assume now that $\lambda:=T^{\prime}(\bar{x},u)\in \Bbb R$. Then
there exist sequences $u_{k}\rightarrow u$ and
$t_{k}\rightarrow0^{+}$ such that $\lambda_{k}:=t_{k}^{-1}\left[
T(\bar{x}+t_{k}u_{k})-T(\bar {x})\right]
\rightarrow\lambda$. In particular, we may assume that
$\lambda_{k}\in \Bbb R$ for $k\geq1$. It follows that
\[
\bar{x}+T(\bar{x})v+t_{k}(u_{k}+\lambda_{k}v)=\bar{x}%
+t_{k}u_{k}+\left[  T(\bar{x})+t_{k}\lambda_{k}\right]  v\in\Omega.
\]
Since $(u_{k}+\lambda_{k}v)\rightarrow u+\lambda v$, we have that
$u+\lambda v\in K(\bar{x}+T(\bar{x})v,\Omega)$, and so
$\varphi_v (u,K(\tilde {x},\Omega))\leq\lambda$. The conclusion
follows.

The proof of the equality in the case where $v\in\Omega_{\infty}$ is left for
the readers. $\hfill\square$\vspace*{0.05in}

A function $\psi:X\rightarrow(-\infty,\infty]$ is called \emph{calm} at
$\bar{x}\in\mbox{\rm dom }\psi$ if there exist constants $\ell\geq0$ and
$\delta>0$ such that
\[
|\psi(x)-\psi(\bar{x})|\leq\ell\Vert x-\bar{x}\Vert\mbox{ for all }x\in
I\!\!B(\bar{x};\delta).
\]
We also say that $\psi$ is \emph{lower calm} at $\bar{x}$ if there exist
constants $\ell\in I\!\!R$ and $\delta>0$ such that
\[
\psi(x)-\psi(\bar{x})\geq\ell\Vert x-\bar{x}\Vert\mbox{ for all }x\in
I\!\!B(\bar{x};\delta).
\]

It is obvious that if $\psi$ is Lipschitz continuous around $\bar{x}$, then it
is (lower) calm at every point in a neighborhood of $\bar{x}$.

\begin{Corollary}\label{c3.9} Let $\bar{x}\in\operatorname*{dom}T
\setminus\Omega$ and let $\tilde{x}:=\Pi_v(\ox;\O)$. Suppose that $T$ is lower calm at
$\ox$. Then
\begin{equation*}
T^{\prime}(\bar{x},u)=\varphi_v(u;K(\tilde {x},\Omega))
\mbox{ for every }u\in X.
\end{equation*}
\end{Corollary}Proof: By definition,
\[
T^{\prime}(\bar{x};u)=\liminf_{t\rightarrow0^{+},\;v\rightarrow u}%
\dfrac{T(\bar{x}+tv)-T(\bar{x})}{t}\geq\liminf_{t\rightarrow0^{+}%
,\;v\rightarrow u}\dfrac{\ell t\Vert v\Vert}{t}=\ell\Vert u\Vert>-\infty,
\]
where $\ell$ is a constant. The equality then follows directly from Lemma
\ref{lm3.8}. $\hfill\square$

\begin{Theorem}\label{t3.11}Let $\bar{x}\in\operatorname*{dom}T
\setminus\Omega$ and let $\tilde{x}:=\Pi_v(\ox;\O)$.
Then
\begin{equation*}
\partial^-T(\ox)\subseteq\{x^*\in X^*\mid  \la x^*, -v\ra=1\}\cap N^{-}
(\tilde{x}; \O).
\end{equation*}
The opposite inclusion
holds if $T$ is lower calm at $\ox$.
\end{Theorem}Proof: Fix any $x^{\ast}\in\partial^{-}T(\bar{x})$. By Lemma
\ref{lm3.8},
\[
\langle x^{\ast},h\rangle\leq T^{\prime}(\bar{x};h)\leq\varphi_{v}%
(h,K(\tilde{x},\Omega))\mbox{ for all }h\in X.
\]
A similar argument to the proof of Theorem \ref{DH1} gives us $\langle
x^{\ast},-v\rangle\leq1$ and $x^{\ast}\in N^{-}(\tilde{x};\Omega)$.

We also have, using Proposition \ref{time property}, that
\[
\langle x^{\ast},v\rangle\leq T^{\prime}(\bar{x},v)\leq\liminf
_{t\rightarrow0^{+}}\dfrac{T(\bar{x}+tv)-T(\bar{x})}{t}=\liminf_{t\rightarrow
0^{+}}\dfrac{T(\bar{x})-t-T(\bar{x})}{t}=-1.
\]
Thus, $\langle x^{\ast},-v\rangle=1$.

Let us prove the opposite inclusion under the lower
calmness; hence $T^{\prime}(\bar{x},h)>-\infty$ for all $h\in
X$. Fix any $x^{\ast}\in X^{\ast}$ such that $\langle
x^{\ast},-v\rangle=1$ and $x^{\ast}\in N^{-}(\tilde{x};\Omega)$. We
will show that
\[
\langle x^{\ast},h\rangle\leq T^{\prime}(\bar{x};h)\mbox{ for all }h\in X.
\]
Fix any $h\in X$. The inequality obviously holds when
$T^{\prime}(\bar {x};h)=+\infty;$ so, assume that
$T^{\prime}(\bar{x};h)\in\mathbb{R}$. By Corollary \ref{c3.9},
\[
h+T^{\prime}(\bar{x};h)v\in K(\tilde{x},\Omega).
\]
It follows that
\[
\langle x^{\ast},h+T^{\prime}(\bar{x};h)v\rangle\leq0.
\]
Thus,
\[
\langle x^{\ast},h\rangle\leq T^{\prime}(\bar{x};h)\langle x^{\ast}%
,-v\rangle.
\]
Therefore, $\langle x^{\ast},h\rangle\leq T^{\prime}(\bar{x};h)$,
and hence $x^{\ast}\in\partial^{-}T(\bar{x})$. $\hfill\square$

The result referring to $\varphi_{v}(\cdot;\Omega)$
which corresponds to Corollary \ref{c3.9} and Theorem \ref{t3.11} is
the following.

\begin{Corollary} Assume that $v\in\Omega_{\infty}$
and $\bar{x}\in X$ is such that
$\varphi(\bar{x})\in\mathbb{R}$, where $\varphi:=
\varphi_{v}(\cdot;\Omega)$. If $\varphi$ is lower calm at
$\bar{x}$, then
\begin{equation*}
\varphi^{\prime}(\bar{x};u)=\varphi_v(u;K(\bar
{x}+\varphi(\bar{x})v,\Omega)) \mbox{ for every }u\in X,
\end{equation*}
and
\[
\partial^-\varphi(\bar{x}) = N^-(\bar
{x}+\varphi(\bar{x})v;\O)\cap\{x^*\in X^*\mid \langle
x^*,-v\rangle=1\}.
\]
\end{Corollary}

\subsection{H\"{o}lder and Fr\'{e}chet sub\-gradients}

Let $\psi:X\rightarrow(-\infty,\infty]$ be an extended real-valued function
and let $\bar{x}\in\mbox{\rm dom }\psi$. Given $s>0$, an element $x^{*}\in
X^{*}$ is called an $s-$H\"{o}lder sub\-gradient of the function $\psi$ at
$\bar{x}$ if there exist $\delta>0$ and $\sigma>0$ such that
\[
\langle x^{*},x-\bar{x}\rangle\leq\psi(x)-\psi(\bar{x})+\sigma\Vert x-\bar
{x}\Vert^{1+s}\mbox{ for all }x\in I\!\!B(\bar{x};\delta).
\]
The set of all $s-$H\"{o}lder sub\-gradients of the function at $\bar{x}$ is
called the $s-$H\"{o}lder subdifferential of the function at this point and is
denoted by $\partial_{s}\psi(\bar{x})$.\vspace*{0.05in}

Similarly, the $s-$H\"{o}lder normal cone to a set $\Omega$ at $\bar{x}%
\in\Omega$ is the set of all $x^{*}\in X^{*}$ such that there exist $\delta>0$
and $\sigma>0$ such that
\[
\langle x^{*},x-\bar{x}\rangle\leq\sigma\Vert x-\bar{x}\Vert^{1+s}\mbox{
for all }x\in I\!\!B(\bar{x};\delta)\cap\Omega.
\]
In the case where $s=1$ and $X$ is a Hilbert space, these structures reduce to the proximal
subdifferential and proximal normal cone, respectively; see \cite{c-1998}.

\begin{Proposition} For any $\ox\in \O$, one has the following representation of
$s-$H${\ddot{o}}$lder sub\-gradients:
\begin{equation}\label{h1}
\partial_s T(\ox)=\{x^*\in X^*\mid  \la x^*,- v\ra \leq 1\} \cap
N_s(\ox;\O).
\end{equation}
Suppose additionally that $v\in\O_\infty$. Then
\begin{equation}\label{h2e}
\partial_s T(\ox)=\{x^*\in X^*\mid  -1\leq \la x^*, v\ra \leq 0\}
\cap N_s(\ox;\O).
\end{equation}
\end{Proposition}Proof: Fix any $x^{*}\in\partial_{s}T(\bar{x})$.
Then there exist $\delta>0$ and $\sigma>0$ such that
\[
\label{fe}\langle x^{*},x-\bar{x}\rangle\leq T(x)-T(\bar
{x})+\sigma\Vert x-\bar{x}\Vert^{1+s}\mbox{ for all }x\in
I\!\!B(\bar{x};\delta).
\]
Since $T(x)=0$ for all $x\in\Omega$, one has
\[
\langle x^{*},x-\bar{x}\rangle\leq\sigma\Vert x-\bar{x}\Vert^{1+s}%
\mbox{ for all }x\in I\!\!B(\bar{x};\delta)\cap\Omega.
\]
This implies $x^{*}\in N_{s}(\bar{x};\Omega)$. Since $\bar{x}-tv\in
I\!\!B(\bar{x};\delta)$ for $t>0$ sufficiently small, one has
\[
\langle x^{*},-tv\rangle\leq T(\bar{x}-tv)+\sigma t^{1+s}\Vert
v\Vert^{1+s}\leq t+\sigma t^{1+s}\Vert v\Vert^{1+s}.
\]
It follows that
\[
\langle x^{*},-v\rangle\leq1+\sigma t^{s}\Vert v\Vert^{1+s}.
\]
Letting $t\rightarrow0$, one has $\langle x^{*},-v\rangle\leq1$. The inclusion
$\subseteq$ in (\ref{h1}) has been proved.

Now fix any $x^{*}\in X^{*}$ such that $\langle x^{*},-v\rangle\leq1$ and
$x^{*}\in N_{s}(\bar{x};\Omega)$. Then there exist $\delta>0$ and $\sigma>0$
such that
\[
\langle x^{*},x-\bar{x}\rangle\leq\sigma\Vert x-\bar{x}\Vert^{1+s}%
\mbox{ for all }x\in I\!\!B(\bar{x};\delta)\cap\Omega.
\]
Suppose by contradiction that $x^{*}\notin\partial_{s}T(\bar{x})$. Then there exist sequences $\sigma_{k}\rightarrow\infty$ and $x_{k}\rightarrow\bar{x}$ such that
\[
\langle x^{*},x_{k}-\bar{x}\rangle>T(x_{k})+\sigma_{k}\Vert
x_{k}-\bar{x}\Vert^{1+s}\geq T(x_{k})=:t_{k}.
\]
This implies $t_{k}\rightarrow0$ as $k\rightarrow\infty$ and $x_{k}\neq\bar
{x}$ for every $k$. Moreover, $t_{k}\leq\Vert x^{*}\Vert\;\Vert x_{k}-\bar
{x}\Vert$. Let $\tilde{x}_{k}:=x_{k}+t_{k}v\in\Omega$. For sufficiently large
$k$, one has
\[
\langle x^{*},\tilde{x}_{k}-\bar{x}\rangle\leq\sigma\Vert\tilde{x}_{k}-\bar
{x}\Vert^{1+s}\leq\sigma(\Vert x_{k}-\bar{x}\Vert+t_{k}\Vert v\Vert)^{1+s}%
\leq\sigma(1+\Vert x^{*}\Vert\;\Vert v\Vert)^{1+s}\Vert x_{k}-\bar{x}%
\Vert^{1+s}.
\]
We also have
\[
\langle x^{*},\tilde{x}_{k}-\bar{x}\rangle=\langle x^{*},x_{k}-\bar{x}%
\rangle+t_{k}\langle x^{*},v\rangle\geq\langle x^{*},x_{k}-\bar{x}%
\rangle-t_{k}>\sigma_{k}\Vert x_{k}-\bar{x}\Vert^{1+s}.
\]
It follows that for sufficiently large $k$,
\[
\sigma_{k}\Vert x_{k}-\bar{x}\Vert^{1+s}<\sigma(1+\Vert x^{*}\Vert\;\Vert
v\Vert)^{1+s}\Vert x_{k}-\bar{x}\Vert^{1+s},
\]
which implies $\sigma_{k}<\sigma(1+\Vert x^{*}\Vert\;\Vert v\Vert)^{1+s}$. We
have arrived at a contradiction since $\sigma_{k}\rightarrow\infty$.

Suppose now that $v\in\Omega_{\infty}$. Using $x:=\bar{x}+tv\in I\!\!B(\bar
{x};\delta)$ for $t>0$ sufficiently small, we arrive at $\langle
x^{*},v\rangle\leq0$, and hence (\ref{h2e}) holds. The proof is now complete.
$\hfill\square$

\begin{Theorem}\label{t3.14} For any $\ox\in \mbox{\rm dom
}T\setminus \O$, one has the following representation of
$s-$H\"{o}lder sub\-gradients:
\begin{equation*}
\partial_s T(\ox)\subseteq\{x^*\in X^*\mid  \la x^*, v\ra = -1\}
\cap N_s(\tilde{x};\O),
\end{equation*}
where $\tilde{x}:=\Pi_v(\ox;\O)$. Suppose further that
$T$ is lower calm at $\ox$. Then
\begin{equation*}
\partial_s T(\ox)=\{x^*\in X^*\mid  \la x^*, v\ra = -1\}
\cap N_s(\tilde{x};\O).
\end{equation*}
\end{Theorem}Proof: Fix any $x^{\ast}\in\partial_{s}T(\bar{x})$. Then
there exist $\delta>0$ and $\sigma>0$ such that
\[
\langle x^{\ast},x-\bar{x}\rangle\leq T(x)-T(\bar{x})+\sigma\Vert
x-\bar{x}\Vert^{1+s}%
\]
for all $x\in I\!\!B(\bar{x};\delta)$. Since $\Omega$ is closed, we can assume
without loss of generality that $I\!\!B(\bar{x};\delta)\cap\Omega=\emptyset$.
Choose $t>0$ sufficiently small so that $\bar{x}-tv\in I\!\!B(\bar{x};\delta
)$. Then (using also Proposition \ref{time property})
\[
\langle x^{\ast},(\bar{x}-tv)-\bar{x}\rangle\leq T(\bar{x}-tv)-T(\bar
{x})+\sigma t^{1+s}\Vert v\Vert^{1+s}\leq T(\bar{x})+t-T(\bar
{x})+\sigma t^{1+s}\Vert v\Vert^{1+s}.
\]
This implies
\[
\langle x^{\ast},-v\rangle\leq1+\sigma t^{s}\Vert v\Vert^{1+s}.
\]
Thus, $\langle x^{\ast},-v\rangle\leq1$. Using $\bar{x}+tv$ in a similar way,
one has $\langle x^{\ast},-v\rangle\geq1$.

Fix any $x\in\Omega$ with $\Vert x-\tilde{x}\Vert<\delta$. Then $\Vert
x-T(\bar{x})v-\bar{x}\Vert<\delta$. Let $t:=T(\bar{x})$. Using
$x-tv+tv\in\Omega$, one has
\begin{align*}
\langle x^{\ast},x-\tilde{x}\rangle &  =\langle x^{\ast},x-tv-\bar{x}\rangle\\
&  \leq T(x-tv)-T(\bar{x})+\sigma\Vert x-\tilde{x}\Vert
^{1+s}\\
&  \leq t-T(\bar{x})+\sigma\Vert x-\tilde{x}\Vert^{1+s}\\
&  =\sigma\Vert x-\tilde{x}\Vert^{1+s}.
\end{align*}
Thus, $x^{\ast}\in N_{s}(\tilde{x};\Omega)$.

Let us prove the opposite inclusion under the calmness of $T$. Let
$\ell,\delta>0$ be such that
\begin{equation}
T(x)-T(\bar{x})\geq-\ell\left\Vert x-\bar{x}\right\Vert \quad\forall
x\in B(\bar{x},\delta). \label{rn1}%
\end{equation}
Take $x^{\ast}\in N_{s}(\tilde{x};\Omega)$ with $\left\langle
x^{\ast },v\right\rangle =-1$. Then there exist $\sigma>0$ and
$\delta^{\prime}>0$ such that
\begin{equation}
\left\langle x^{\ast},u-\tilde{x}\right\rangle \leq\sigma\left\Vert
u-\tilde{x}\right\Vert ^{1+s}\quad\forall u\in B(\tilde{x}%
,\delta^{\prime}). \label{rn3}%
\end{equation}
Assume that $x^{\ast}\notin\partial_{s}T(\bar{x})$. Then there
exists a
sequence $(x_{k})\subseteq B(\bar{x},\delta)$ such that $x_{k}%
\rightarrow\bar{x}$ and
\begin{equation}
\left\langle x^{\ast},x_{k}-\bar{x}\right\rangle >T(x_{k})-T(\bar
{x})+k\left\Vert x_{k}-\bar{x}\right\Vert ^{1+s}\quad\forall k\geq1.
\label{rn2}%
\end{equation}
From (\ref{rn1}) and (\ref{rn2}) we obtain that $\left\vert T(x_{k}%
)-T(\bar{x})\right\vert \leq\ell^{\prime}\left\Vert
x_{k}-\bar {x}\right\Vert $, where
$\ell^{\prime}:=\max(\ell,\left\Vert x^{\ast }\right\Vert )$, and so
$T(x_{k})\rightarrow T(\bar{x})$. Set
$u_{k}:=x_{k}+T(x_{k})v\in\Omega;$ clearly,
$u_{k}\rightarrow\tilde{x}$, and so $u_{k}\in
B(\tilde{x},\delta^{\prime})$ for $k\geq k_{0}$ with
$k_{0}\geq1$ fixed. It follows that for $k\geq k_{0}$ we have
\begin{align*}
\left\langle x^{\ast},x_{k}-\bar{x}\right\rangle -\left[  T(x_{k}%
)-T(\bar{x})\right]   &  =\left\langle x^{\ast},u_{k}-\tilde
{x}\right\rangle \leq\sigma\left\Vert u_{k}-\tilde{x}\right\Vert ^{1+s}\\
&  =\sigma\left\Vert x_{k}-\bar{x}+\left[  T(x_{k})-T(\bar
{x})\right]  v\right\Vert ^{1+s}\\
&  \leq\sigma(1+\ell^{\prime}\left\Vert v\right\Vert )^{1+s}\left\Vert
x_{k}-\bar{x}\right\Vert ^{1+s}.
\end{align*}
Using (\ref{rn2}) we get
\[
k\left\Vert x_{k}-\bar{x}\right\Vert ^{1+s}<\left\langle x^{\ast}%
,x_{k}-\bar{x}\right\rangle -\left[  T(x_{k})-T(\bar{x})\right]
\leq\sigma(1+\ell^{\prime}\left\Vert v\right\Vert )^{1+s}\left\Vert
x_{k}-\bar{x}\right\Vert ^{1+s}%
\]
for $k\geq k_{0}$, whence the contradiction $k<\sigma(1+\ell^{\prime
}\left\Vert v\right\Vert )^{1+s}$ for every $k\geq k_{0}$. Hence $x^{\ast}%
\in\partial_{s}T(\bar{x})$. $\hfill\square$

Similar proofs yield the representations below for Fr\'{e}chet sub\-gradients
of the directional minimal time function (\ref{s}) in both in-set case and
out-of-set case.

\begin{Proposition}\label{fr1}
For any $\ox\in \O$, one has the following representation of Fr\'{e}chet
sub\-gradients:
\begin{equation}\label{f1}
\Hat\partial T(\ox)=\{x^*\in X^*\mid  \la x^*,- v\ra \leq 1\}
\cap \Hat N(\ox;\O).
\end{equation}
Suppose additionally that $v\in\O_\infty$. Then
\begin{equation}\label{f2e}
\Hat\partial T(\ox)=\{x^*\in X^*\mid  -1\leq \la x^*, v\ra
\leq 0\} \cap \Hat N(\ox;\O).
\end{equation}
\end{Proposition}

\begin{Theorem}\label{t3.16} For any $\ox\in \mbox{\rm dom }
T\setminus \O$, one has the
following representation of Fr\'{e}chet sub\-gradients:
\begin{equation*}
\Hat\partial T(\ox)\subseteq\{x^*\in X^*\mid  \la x^*, v\ra =
-1\} \cap \Hat N(\tilde{x};\O),
\end{equation*}
where $\tilde{x}:=\Pi_v(\ox;\O)$. Suppose further that
$T_v(\cdot;\O)$ is lower calm at $\ox$. Then
\begin{equation*}
\Hat\partial T(\ox)=\{x^*\in X^*\mid  \la x^*, v\ra = -1\}
\cap \Hat N(\tilde{x};\O).
\end{equation*}
\end{Theorem}

The result referring to $\varphi_{v}(\cdot;\Omega)$
which corresponds to  Theorems \ref{t3.14} and \ref{t3.16} is the
following.

\begin{Corollary} Assume that $v\in\Omega_{\infty}$,
$\bar{x}\in X$ is such that
$\varphi(\bar{x})\in\mathbb{R}$, where $\varphi:=
\varphi_{v}(\cdot;\Omega)$, and $s>0$. Then
\begin{equation*}
\partial_s \varphi(\bar{x})\subseteq\{x^*\in X^*\mid  \la x^*, v\ra =
-1\} \cap  N_s(\bar{x}+\varphi(\bar{x})v;\O),
\end{equation*}
and
\begin{equation*}
\Hat\partial \varphi(\bar{x})\subseteq\{x^*\in X^*\mid  \la
x^*, v\ra = -1\} \cap \Hat
N(\bar{x}+\varphi(\bar{x})v;\O).
\end{equation*}
Moreover, if $\varphi$ is lower calm at $\bar{x}$, then
equalities hold in the previous two inclusions.
\end{Corollary}

\subsection{Limiting Sub\-gradients}

\begin{Theorem}\label{limiting1} For any $\ox\in \O$, one has the
following representation of
limiting sub\-gradients:
\begin{equation*}
\partial T(\ox)=\{x^*\in X^*\mid  \la x^*, -v\ra \leq 1\} \cap N(\ox;\O).
\end{equation*}
\end{Theorem}Proof: Fix any $x^{*}\in\partial T(\bar{x})$. Then
there exist sequences $x_{k}^{*}\xrightarrow{w^*}x^{*}$ and $x_{k}%
\xrightarrow{T}\bar{x}$ with $x_{k}^{*}\in\widehat{\partial}T(x_{k})$. Let
$\tilde{x}_{k}:=\Pi_v(x_{k};\Omega)$. Then $\tilde{x}_{k}\xrightarrow{\O}\bar
{x}$ and $x_{k}^{*}\in\widehat{N}(\tilde{x}_{k};\Omega)$. This implies
$x^{*}\in N(\bar{x};\Omega)$. In both cases: $x_{k}\in\Omega$ and $x_{k}%
\notin\Omega$, we always have $\langle x_{k}^{*},-v\rangle\leq1$. Thus
$\langle x^{*},-v\rangle\leq1$.

Let us prove the opposite inclusion. Fix $x^{*}\in X^{*}$ with $\langle
x^{*},-v\rangle\leq1$ and $x^{*}\in N(\bar{x};\Omega)$. Then there exist
sequences $x_{k}^{*}\xrightarrow{w^*}x^{*}$ and $x_{k}\xrightarrow{\O}\bar{x}$
with $x_{k}^{*}\in\widehat{N}(x_{k};\Omega)$. Define $\gamma_{k}:=\langle
x_{k}^{*},-v\rangle$. If $\gamma_{k}\leq1$ for a subsequence (without
relabeling), then $x_{k}^{*}\in\widehat{\partial}T(x_{k})$, and
hence $x^{*}\in\partial T(\bar{x})$. So we can assume $\gamma
_{k}>1$ for every $k$. Clearly, $\gamma_{k}\rightarrow\gamma:=\langle
x^{*},-v\rangle=1$. Let $\tilde{x}_{k}^{*}:=\dfrac{x_{k}^{*}}{\gamma_{k}}$.
Then $\tilde{x}_{k}^{*}\in\widehat{\partial}T(x_{k})$. So again,
$x^{*}\in\partial T(\bar{x})$.$\hfill\square$

We say that $\Omega$ satisfies \emph{property $P$ around $\bar{x}$} with a radius $r>0$
if there exists a neighborhood $V$ of $\bar{x}$ such that $x-tv\notin\Omega$
for all $t\in(0,r]$ and for all $x\in V\cap\mbox{\rm bd }\Omega$.

For example, if $\Omega$ is the epigraph of a continuous function
$\psi:X\rightarrow(-\infty,\infty]$ and $v=(0,1)$, where $0$ is the zero
element of $X$, then property $P$ is satisfied.

\begin{Theorem}\label{limiting 2} For any $\ox\in
\mbox{\rm dom }T\setminus \O$, one has the following upper estimate of
limiting sub\-gradients:
\begin{equation*}
\partial T(\ox)\subseteq\{x^*\in X^*\mid  \la x^*, v\ra = -1\} \cap
N(\tilde{x};\O),
\end{equation*}
where $\tilde{x}:=\Pi_v(\ox; \O)$.
Suppose further that $T$ is lower calm around $\ox$. Then
\begin{equation*}
\partial T(\ox)=\{x^*\in X^*\mid  \la x^*, v\ra = -1\} \cap
N(\tilde{x};\O),
\end{equation*}
under the assumption that $\O$ satisfies condition \mbox{\rm $P$}
around $\tilde{x}$ with the radius $r=T(\ox)$.
\end{Theorem}Proof: Fix any $x^{\ast}\in\partial T(\bar{x})$. Then
there exist sequences $x_{k}\rightarrow\bar{x}$ with $T(x_{k})\rightarrow T(\bar{x})$ and $x_{k}^{\ast}%
\xrightarrow{w^*}x^{\ast}$ with $x_{k}^{\ast}\in\widehat{\partial}T(x_{k})$. Under the assumption made, $\langle x_{k}^{\ast},v\rangle=-1$
and $x_{k}^{\ast}\in\widehat{N}(\tilde{x}_{k};\Omega)$ for sufficiently large
$k$, where $\tilde{x}_{k}:=x_{k}+T(x_{k})v\in\Omega$ (see Theorem
\ref{t3.16}). Then $\langle x^{\ast},v\rangle=-1$. Clearly, $\tilde{x}%
_{k}\overset{\Omega}{\rightarrow}\tilde{x}$. Thus, $x^{\ast}\in\widehat
{N}(\tilde{x};\Omega)$. The first inclusion has been proved.

Let us prove the opposite inclusion. Fix any $x^{\ast}\in X^{\ast}$ with
$\langle x^{\ast},v\rangle=-1$ and $x^{\ast}\in\widehat{N}(\tilde{x};\Omega)$.
Then there exist $\tilde{x}_{k}\xrightarrow{\O}\tilde{x}$ and $x_{k}^{\ast}%
\in\widehat{N}(\tilde{x}_{k};\Omega)$ with $x_{k}^{\ast}%
\xrightarrow{w^*}x^{\ast}$. Let $\gamma_{k}:=-\langle x_{k}^{\ast},v\rangle$.
Then $\gamma_{k}\rightarrow1$ as $k\rightarrow\infty$. For sufficiently large
$k$, let $\tilde{x}_{k}^{\ast}:=\dfrac{x_{k}^{\ast}}{\gamma_{k}}$. Then
$\langle\tilde{x}_{k}^{\ast},v\rangle=-1$ and $\tilde{x}_{k}^{\ast}\in
\widehat{N}(\tilde{x}_{k};\Omega)$ by the cone property of the Fr\'{e}chet
normal cone. Clearly, $\tilde{x}_{k}$ belongs to the boundary of $\Omega$. Let
$x_{k}:=\tilde{x}_{k}-rv$. Under the P
property, and using Proposition \ref{time property}, we have that $T(x_{k})=r$
and $\Pi(x_{k};\Omega)=\tilde{x}_{k}$. Thus, by Theorem \ref{t3.16},
$\tilde{x}_{k}^{\ast}\in\widehat{\partial}T(x_{k})$. Since
$x_{k}\xrightarrow{T}\bar{x}$, one has $x^{\ast}\in\partial
T(\bar{x})$. $\hfill\square$

\begin{Theorem}\label{singular 1}For $\ox\in\O$, one has
\begin{equation*}
\partial^\infty T(\ox)= \{v\}^+\cap N(\ox;\O).
\end{equation*}
\end{Theorem}Proof: The proof of $\subseteq$ follows directly from the
definition of singular sub\-gradients and the proof of Theorem
\ref{limiting 2}. Let us prove the opposite inclusion.
Fix any $x^{\ast}\in\{v\}^{+}\cap
N(\bar{x};\Omega)$. Then there exist $x_{k}\xrightarrow{\O}\bar{x}$,
$x_{k}^{\ast}\xrightarrow{w^*}x^{\ast}$ such that
$x_{k}^{\ast}\in\widehat{N}(x_{k};\Omega)$; clearly
$x_{k}\xrightarrow{T}\bar {x}$. Set $\mu_{k}:=\left\langle
x_{k}^{\ast},v\right\rangle \rightarrow \left\langle
x^{\ast},v\right\rangle \geq0$. If $\mu_{k}\geq0$ for a subsequence
(without relabeling), then $\langle kx_{k}^{\ast},v\rangle
\geq0\geq-1$, and hence
$kx_{k}^{\ast}\in\widehat{\partial}T(x_{k})$. This
implies $x^{\ast}\in\partial^{\infty}T(\bar{x})$. In the
contrary case $\mu_{k}<0$ for every $k$, and so $\lambda_{k}:=-\mu
_{k}\downarrow0$. Then
$\langle\dfrac{x_{k}^{\ast}}{\lambda_{k}},v\rangle=-1$, and so
$\dfrac{x_{k}^{\ast}}{\lambda_{k}}\in\widehat{\partial}T(x_{k})$.
Then, by definition, again
$x^{\ast}\in\partial^{\infty}T(\bar{x})$.
$\hfill\square$ \vspace*{0.05in}

The proof of the theorem below is also straightforward.

\begin{Theorem}\label{soutset} For any $\ox\in
\mbox{\rm dom }T\setminus \O$, one has
\begin{equation*}
\partial^\infty T(\ox) \subseteq\{v\}^\perp\cap N(\ox+T(\ox)v;\O).
\end{equation*}
\end{Theorem}

\section{Lipschitz properties of directional minimal time functions}

In this section, we are going to study Lipschitz properties of the directional
minimal time function (\ref{s}). Necessary and sufficient conditions for
globally Lipschitz property and locally Lipschitz property/Lipschitz
continuity will be established.

\begin{Proposition}\label{l1} Suppose that $v\in\mbox{\rm int }\O_\infty$.
Then $T$ is globally Lipschitz with Lipschitz constant
$$\ell:=\inf\{r^{-1}\mid r>0,\ \B(v; r)\subseteq \O_\infty \}=\frac{1}{\mbox{\rm dist}(v, \mbox{\rm bd }\O_\infty)}.$$
\end{Proposition}Proof: Let $K:=\Omega_{\infty}$. We first show that
$\mbox{\rm dom }T_{v}(\cdot;\Omega)=X$. Indeed, since $v\in\mbox{\rm int }K$,
\[
\mbox{\rm dom }T_{v}(\cdot;\Omega)=\Omega-\mbox{\rm span }\{v\}=\Omega
+[K-\mbox{\rm span }\{v\}]=\Omega+X=X.
\]
It follows from Proposition \ref{ll1} (2) that
\[
T_{v}(x;\Omega+K)\leq T_{v}(y;\Omega)+T_{v}(x-y;K).
\]
Then
\[
T_{v}(x;\Omega)\leq T_{v}(y;\Omega)+T_{v}(x-y;K).
\]
This implies
\begin{equation}
T_{v}(x;\Omega)-T_{v}(y;\Omega)\leq T_{v}(x-y;K). \label{lip11}%
\end{equation}
Fix any $r>0$ such that $I\!\!B(v;r)\subseteq K$. Then $I\!\!B(0;r)\subseteq
K-v$. We have the following for $u\notin K$ by the cone property of $K$
\begin{align*}
T_{v}(u;K)  &  =\inf\{t\geq0\mid u+tv\in K\}\\
&  =\inf\{t>0\mid u\in t(K-v)\}\\
&  \leq\inf\{t>0\mid u\in tI\!\!B(0;r)\}=\Vert r^{-1}u\Vert.
\end{align*}
Since $T_v(u; K)=0\leq r^{-1}\Vert u\Vert$ for $u\in K$, one has
\begin{equation*}
T_v(u; K)\leq r^{-1}\Vert u\Vert \quad \forall u\in X.
\end{equation*}
It follows from (\ref{lip11}) that
\[
T_{v}(x;\Omega)-T_{v}(y;\Omega)\leq T_{v}(x-y;K)\leq\ell\Vert x-y\Vert.
\]
That implies
\[
|T_{v}(y;\Omega)-T_{v}(x;\Omega)|\leq\ell\Vert x-y\Vert.
\]
The proof is now complete. $\hfill\square$

\begin{Proposition}  The function
$T$ is finite-valued and
Lipschitz if and only if $v\in \mbox{\rm int }\O_\infty$.
\end{Proposition}
Proof: Suppose $T$ is
finite-valued and Lipschitz. Then $v\in\Omega_{\infty}$. Otherwise,
there exist $\bar{x}\in\Omega$ and $t>0$ such that
$\bar{x}+tv\notin\Omega$. Because $\bar{x}+tv\in
X=\operatorname*{dom}T=\Omega-\mathbb{R}_{+}v$, there exists $s>0$
such that $\bar{x}+tv+sv\in\Omega$. Taking $\bar {t}:=t+s$
and using Corollary \ref{cor2.9}, we get some $t_{0}\in
\lbrack0,\bar{t})$ such that $T$ is not continuous at
$\bar{x}+t_{0}v$, contradicting our hypothesis. Therefore,
$v\in\Omega_{\infty}.$

Let $\ell\geq0$ be the Lipschitz constant of $T$, that is
\[
|T(x)-T(y)|\leq\ell\Vert u-v\Vert\mbox{ for all }x,y\in X.
\]
We are going to show that $v+\dfrac{1}{\ell}I\!\!B\in\Omega_{\infty}$. Indeed,
fix an $e\in\dfrac{1}{\ell}I\!\!B$. Using Proposition \ref{time property}, for
any $x\in\Omega$ and $t\geq0$, one has
\[
T(x+t(v+e))=T(x+te+tv)=\max\{T(x+te)-t,0\}\leq\max\{T(x)+t\ell\Vert
e\Vert-t,0\}\leq0.
\]
This implies $x+t(v+e)\in\Omega$. Thus, $v+e\in\Omega_{\infty}$.

The converse follows from Proposition \ref{l1}. $\hfill\square$ \vspace
*{0.05in}

In what follows, we are going to characterize the Lipschitz continuity of the
minimal time function (\ref{s}) using both direct and generalized
differentiation approaches.

\begin{Lemma}\label{lip1} Suppose that $T$ is Lipschitz
continuous around $\tilde{x}:=\Pi_v(\ox;\O)$, where $\ox\in
\operatorname*{dom}T$. Then $T$ is
Lipschitz continuous around $\ox$. Moreover,  the
converse holds if $v\in\Omega_\infty$.
\end{Lemma}
Proof: By hypothesis, there exist $\ell\geq0,$
$\delta>0$ such that $T$ is finite on
$I\!\!B(\tilde{x};\delta)$ and
\[
|T(x)-T(y)|\leq\ell\Vert x-y\Vert\mbox{ for all }x,y\in
I\!\!B(\tilde{x};\delta).
\]
If $\ox\in \O$ the conclusion is obvious. In the contrary case take
$0<\varepsilon<\min\{T(\bar{x}),\delta/(1+\left\Vert
v\right\Vert )\}$. Then there exists
$\delta^{\prime}\in(0,\varepsilon]$ such that
$T(x)>t:=T(\bar{x})-\varepsilon>0$ for every $x\in
B(\bar{x};\delta^{\prime})$. Then for $x\in
B(\bar{x};\delta^{\prime})$ we have that
\[
\left\Vert (x+tv)-\tilde{x}\right\Vert =\left\Vert (x-\bar
{x})-\varepsilon v\right\Vert \leq\varepsilon(1+\left\Vert v\right\Vert
)<\delta,
\]
and so $x+tv\in I\!\!B(\tilde{x};\delta)$, and $T(x+tv)=T(x)-t$ (by
Proposition \ref{time property}). It follows that
\[
\left\vert T(x)-T(y)\right\vert =\left\vert T(x+tv)+t-T(y+tv)-t\right\vert
=\left\vert T(x+tv)-T(y+tv)\right\vert \leq\ell\left\Vert x-y\right\Vert
\]
for $x,y\in B(\bar{x};\delta^{\prime})$.

Assume that $v\in\Omega_{\infty}$ and that $T$ is
Lipschitz around $\bar{x}$. Let $\bar{x}\notin\Omega.$
Using (\ref{t-phi}) (see also Remark \ref{rem2.5}), we obtain that
$\varphi$ is Lipschitz around $\bar{x}$. Using now
(\ref{p-phi}) we obtain that $\varphi$ is Lipschitz around
$\tilde{x}$. Using again (\ref{t-phi}) we get the Lipschitz
continuity of $T.$ The proof is complete. $\hfill\square$

\begin{Lemma}\label{lip2} Let $\ox\in\mbox{\rm dom }T$. The function $T$ is
Lipschitz continuous around $\ox$ if and only if there exist
$\ell\geq 0$ and $\delta >0$ such that
\begin{equation}
T(x+u)\leq \ell \|u\| \mbox{ for all }u\in \delta \B \mbox{ and
}  x\in \O\cap \B(\ox; \delta ). \label{rlip2}
\end{equation}
\end{Lemma}
Proof: The implication $\mathbb{\Rightarrow}$ is
obvious. Assume (\ref{rlip2}) holds for $\ell\geq0$ and $\delta>0.$
Consider $\delta^{\prime }:=\delta/(2+\ell\left\Vert v\right\Vert )$
and take $u,u^{\prime}\in \delta^{\prime}I\!\!B$. Then $\left\Vert
u-u^{\prime}\right\Vert \leq 2\delta^{\prime}\leq\delta$ and
$T(\bar{x}+u) \leq\ell\delta^{\prime}$. We may assume that
$(\ell\delta^{\prime}\geq)$ $T(\bar{x}+u)\geq t:=T(\bar
{x}+u^{\prime})$ $(\geq0)$. Using Proposition \ref{time property},
we have that $T(\bar{x}+u+tv)=T(\bar{x}+u)-t$. Since
$\left\Vert u+tv\right\Vert \leq\delta^{\prime}+t\left\Vert
v\right\Vert \leq \delta^{\prime}\left(  1+\ell\left\Vert
v\right\Vert \right)  \leq\delta$ and
$\bar{x}+u^{\prime}+tv\in\Omega\cap
I\!\!B(\bar{x};\delta)$, using (\ref{rlip2}) we get\[
\left\vert T(\bar{x}+u)-T(\bar{x}+u^{\prime})\right\vert
=T(\bar{x}+u)-t=T(\bar{x}+u+tv)
\leq\ell\left\Vert u-u^{\prime}\right\Vert .
\]
Therefore, $T$ is Lipschitz on
$I\!\!B(\bar{x};\delta^{\prime})$ with the same constant
$\ell$. \hfill$\square$ \vspace*{0.05in}

Recall that $\Omega$ is epi-Lipschitz at $\bar{x}\in X$ in the direction
$v\neq0$ if there exists $\delta>0$ such that for all $\omega\in\Omega\cap
I\!\!B(\bar{x};\delta)$, $u\in I\!\!B(v;\delta)$, and $\lambda\in
\lbrack0,\delta]$, one has $\omega+\lambda u\in\Omega$; see, e.g.,
\cite{rw,tz}.

\begin{Theorem}\label{lip3}Let $\ox\in \O$. If $\O$ is epi-Lipschitz
at $\ox$ in the direction $v$, then $T$ is Lipschitz
continuous around $\ox$. Moreover,  the converse
holds if $v\in\Omega_\infty$.
\end{Theorem}Proof: We only need to show that (\ref{rlip2}) holds for some
$\ell\geq0$ and $\delta>0$. Under the epi-Lipschitz condition, we see that
$T$ is finite around $\bar{x}$. By contradiction, there
exist $v_{k}\rightarrow0$, $x_{k}\xrightarrow{\O}\bar{x}$ and
\begin{equation}
T(x_{k}+v_{k})>k\Vert v_{k}\Vert. \label{ct}%
\end{equation}
Then $v_{k}\neq0$ for every $k$. With the same notation in the definition of
the epi-Lipschitz property, one has that
\[
x_{k}+t(v+e)\in\Omega\mbox{ for }t\in\lbrack0,\delta],\Vert e\Vert\leq
\delta\mbox{ and sufficiently large }k.
\]
This implies $T(x_{k}+te)\leq t$. Thus,
\[
T(x_{k}+v_{k})=T(x_{k}+\dfrac{\Vert v_{k}\Vert}{\delta}%
\dfrac{\delta v_{k}}{\Vert v_{k}\Vert};\Omega)\leq\dfrac{\Vert v_{k}\Vert
}{\delta}\mbox{ for sufficiently large }k.
\]
Comparing with (\ref{ct}), we have arrived at a contradiction.

Suppose that $T$ is Lipschitz continuous around $\bar{x}%
\in\mbox{\rm dom }T$ and $v\in\Omega_{\infty}$. Suppose by
contradiction that $\Omega$ is not epi-Lipschitz around $\bar{x}$ in the
direction $v$. Then there exist sequences $x_{k}\xrightarrow{\O}\bar{x}$,
$u_{k}\rightarrow0$, $t_{k}\rightarrow0$, $t_{k}>0$, such that $x_{k}%
+t_{k}(u_{k}+v)\notin\Omega$. Then $t_{k}<T(x_{k}+t_{k}u_{k})$
because $v\in
\O_\infty$. So
\[
t_{k}<T(x_{k}+t_{k}u_{k})=T(x_{k}+t_{k}u_{k})-T(x_{k})\leq\ell t_{k}\Vert u_{k}\Vert.
\]
Since $u_{k}\rightarrow0$, we have arrived at a contradiction. $\hfill\square$

\begin{Corollary}\label{lip4}Let $\ox\notin \O$ such that $T(\ox)$
is finite and let $\tilde{x}:=\Pi(\ox;\O)$. If $\O$ is epi-Lipschitz
at $\tilde{x}$ in the direction $v$, then $T$ is
Lipschitz continuous around $\ox$. Moreover,  the
converse holds if $v\in\Omega_\infty$.
\end{Corollary}
Proof: The implication $\mathbb{\Rightarrow}$ follows directly from Theorem
\ref{lip3} and Lemma \ref{lip1}.

Let $v\in\Omega_{\infty}$ and assume that $T$ is
Lipschitz continuous around $\bar{x}$. By Lemma \ref{lip1} we
have that $T$ is Lipschitz continuous around
$\tilde{x}\in\Omega$. By Theorem \ref{lip3} we obtain that
$\Omega$ is epi-Lipschitz at $\tilde{x}$ in the direction $v.$
$\hfill\square$\vspace*{0.05in}

In \cite[Theorem 7]{tz}, a necessary and sufficient condition for the
Lipschitz continuity of the scalarization function was proved under the
free-disposal condition, that is $\O+P=\O$ for some closed convex cone $P$ with $v\in P$. Notice that in Theorem \ref{lip3} and Corollary
\ref{lip4}, the free-disposal condition is not required to prove the
sufficient condition for Lipschitz continuity of the directional minimal time
function (\ref{s}). The proof of \cite[Theorem 7]{tz} is not applicable to our
results since it is based on time property from \cite[Theorem 2.3.1]{z}, which
is not satisfied by the directional minimal time function. \vspace*{0.05in}

In the theorem below, we are able to fully characterize the Lipschitz
continuity of the directional minimal time function (\ref{s}) without using
the free-disposal condition. For simplicity, we present our results in finite
dimensions.
%\vspace*{0.05in}

\begin{Theorem}\label{ls1} Let $X$ be a finite dimensional space and let
$\ox\in\operatorname*{dom}T$.

$(1)$ Assume that $\ox \in \O$; then $T$ is Lipschitz
continuous around $\ox$ if and only if $\{v\}^+\cap
N(\ox;\O)=\{0\}$.

$(2)$ Assume that $\ox \notin \O$ and set $\tilde x:=\Pi(\ox;\O)$.
If $\{v\}^+\cap N(\tilde{x};\O)=\{0\}$, then $T$ is
Lipschitz continuous around $\ox$. The converse holds true if $v\in
\O_\infty$.

\end{Theorem}
Proof: $(1)$ On one hand we have
$\partial^{\infty}T(\bar{x})=\{v\}^{+}\cap
N(\bar{x};\Omega)$ by Theorem \ref{singular 1}. On the other hand,
by \cite[Theorem 9.13]{rw}, we have that $T$ is
Lipschitz continuous around $\bar{x}$ iff
$\partial^{\infty}T(\bar {x})=\{0\}$. The conclusion (i)
follows.

$(2)$ In this case we have
$\partial^{\infty}T(\bar{x})\subseteq \{v\}^{+}\cap
N(\tilde{x};\Omega)$ by Theorem \ref{soutset}, and so $T$ is Lipschitz
continuous around $\bar{x}$ as in (i). Assume that $v\in \O_\infty$
and $T$ is Lipschitz continuous around $\ox$. By Lemma
\ref{lip1}, we have that $T$ is Lipschitz continuous
around $\tilde{x}\in \O$. By $(1)$ we obtain that $\{v\}^\perp\cap
N(\tilde{x};\O)=\{0\}$.  $\hfill\square$

\begin{Example}{\rm (1) Let $\O:=\{(x,y)\in {\Bbb R}^2\mid  y\geq -|x|\}$
and let $v=(0,1)$. Then $T$ is Lipschitz continuous at $\ox=(0,0)\in
\O$ since the condition $\{v\}^{+}\cap
N(\bar{x};\Omega)=\{0\}$ is satisfied.\\[1ex]
(2) Let $\O:=\{(x,y)\in {\Bbb R}^2\mid  y\geq -|x|^{1/2}\}$ and let
$v=(0,1)$. Then $T$ is not Lipschitz continuous at $\ox=(0,0)\in \O$
since the condition $\{v\}^{+}\cap N(\bar{x};\Omega)=\{0\}$ is
violated.\\[1ex]
(3) The converse of Theorem \ref{ls1}~$(2)$ does not hold true in
general. Indeed, let $v=(0, 1)$ and
\begin{equation*}
\O:=\{(x,y)\in \Bbb R^2\mid  y= 0\}\cup\{(x,y)\in \Bbb R^2\mid x=0,
y\geq 0\}.
\end{equation*}
Then $T$ is Lipschitz continuous around
$\ox=(0,-2)\notin \O$. However, $$\{v\}^\perp \cap
N(\ox;\O)=\{(x,y)\in \Bbb R^2\mid y=0\}\neq \{0\}.$$}
\end{Example}

\section{Applications to location problems}

In this section, we are going to apply the results obtained previously to
study directional location problems. To the best of our knowledge, the
location model of this type has not been considered in the literature.
\vspace*{0.05in}

Given the nonempty closed target sets $\Omega_{i}$ for $i=1,\ldots,n$ and $n$
directions $v_{i}\neq0$ for $i=1,\ldots,n$, and given a nonempty closed
constraint set $\Omega_{0}$, find a point $\bar{x}\in\Omega_{0}$ such that the
sum of the times to reach the target sets is minimal. The optimization model
is
\begin{equation}
\mbox{\rm minimize }S(x):=\sum_{i=1}^{n}T_{v_{i}}(x;\Omega_{i}%
)\mbox{ subject to }x\in\Omega_{0}. \label{p}%
\end{equation}
It is clear that
\[
\mbox{\rm dom }S=\cap_{i=1}^{n}\mbox{\rm dom }T_{v_{i}}(\cdot;\Omega_{i}%
)=\cap_{i=1}^{n}[\Omega_{i}-\mbox{\rm cone}\,\{v_{i}\}].
\]
\begin{Proposition} Suppose $\mbox{\rm dom }S\cap \O_0\neq \emptyset$.
Then the optimization problem {\rm(\ref{p})} has an optimal solution
under one of the following conditions:\\
{\rm (1)} At least one of the sets among $\O_i$ for $i=0,1,\ldots,n$
is compact.\\
{\rm (2)} $X$ is a reflexive Banach space, $\O_i$ for
$i=0,1,\ldots,n$ are convex and at least one of them is bounded.
\end{Proposition}
Proof: Under the assumptions made, one has
\[
\gamma:=\inf\{S(x)\mid x\in\Omega_0\}<\infty.
\]
Let us first suppose that (1) is satisfied. In the case where the
constraint $\Omega_0$ is compact, an optimal solution exists by the
classical Weierstrass theorem since $S$ is lower semicontinuous.
Suppose without loss of generality that $\Omega_{1}$ is compact. Let
$(x_{k})\subseteq\Omega_0$ be a minimizing sequence. That means
$S(x_{k})\rightarrow\gamma$ as $k\rightarrow\infty$. Thus,
$T_{v_{1}}(x_{k};\Omega_{1})<\gamma+1$ for sufficiently large $k$.
For $t_{k}:=T_{v_{1}}(x_{k};\Omega_{1})$, one has
\[
x_{k}+t_{k}v_{1}\in\Omega_{1}.
\]
Since $\Omega_{1}$ is compact, it is clear that $(x_{k})$ has a convergent
subsequence (without relabeling) to $\bar{x}\in\Omega_0$. Since $S$ is lower
semicontinuous,
\[
S(\bar{x})\leq\liminf_{k\rightarrow\infty}S(x_{k})=\gamma.
\]
Therefore, $\bar{x}$ is an optimal solution of the problem.

In the case where (2) is satisfied, we use a similar argument using the
observation that $S$ is weakly lower semicontinuous since it is convex and
lower semicontinuous. Moreover, every closed bounded convex set in a reflexive
Banach space is weakly sequentially compact. $\hfill\square$

\begin{Theorem} Suppose that $\O_0$ is convex and that $\O_i$ are strictly
convex for $i=1,\ldots,n$, $n\geq 2$. Suppose that any interval $[x,y]$,
$x,y\in\O_0$, $x\neq y$, does not intersect at least two sets among
$\O_i$ for $i=1,\ldots,n$. Suppose further that any set of two
vectors $\{v_i, v_j\}$, $i\neq j$, is linearly independent. Then the
optimization problem {\rm (\ref{p})} has at most one optimal
solution.
\end{Theorem}
Proof: We will show that $S$ is strictly convex on $\Omega
_{0}\cap\mbox{\rm dom }S$. Suppose by contradiction that there exist
$a,b\in\Omega_0\cap\mbox{\rm dom }S$, $a\neq b$, and $t\in(0,1)$
such that
\[
S(ta+(1-t)b)=tS(a)+(1-t)S(b).
\]
This implies
\[
T_{v_{i}}(ta+(1-t)b)=tT_{v_{i}}(a)+(1-t)T_{v_{i}}(b)\mbox{ for all }i=1,\ldots
,n.
\]
Suppose $[a,b]\cap\Omega_{i}=\emptyset$ and $[a,b]\cap\Omega_{j}=\emptyset$
for $i\neq j$. Then $b-a\in\mbox{\rm span }\{v_{i}\}\cap
\mbox{\rm span }\{v_{i}\}=\{0\}$, so $a=b$ by Proposition \ref{2.5}, which is
a contradiction.$\hfill\square$ \vspace*{0.05in}

In what follows we are going to establish necessary and sufficient optimality
condition for problem (\ref{p}). For every $u\in X$, define
\[
I(u):=\{i\in\{1,\ldots,n\}\mid u\in\Omega_{i}\},
\]
and
\[
J(u):=\{i\in\{1,\ldots,n\}\mid u\notin\Omega_{i}\}.
\]
\begin{Theorem}\label{op1} Let $X$ be an Asplund space (see \cite{mor} for the definition). Consider the
optimization problem {\rm (\ref{p})}. Suppose that $\ox\in \O_0$ is
an optimal solution of the problem and $T_{v_i}(\cdot; \O_i)$ is
Lipschitz continuous around $\ox$ for all $i=1,\ldots,n$. Then there
exist $x^*_i\in
X^*$ with the following properties:\\
{\rm (1)} $x^*_i\in N(\tilde{x}_i; \O_i)$, where $\tilde{x}_i:=\Pi(\ox;\O_i)$
for $i=1,\ldots,n.$\\
{\rm (2)} $\la x_i^*, -v_i\ra \leq 1$ for all $i\in I(\ox)$ and $\la x_j^*,
v_j\ra =-1$ for all $j\in J(\ox)$.\\
{\rm (3)} $-\sum_{i=1}^n x^*_i\in N(\ox; \O_0)$.
\end{Theorem}
Proof: It is clear that $\bar{x}\in\Omega_0$ is an optimal solution
of the optimization problem (\ref{p}) if and only if it is a
solution of the following unconstrained optimization problem:
\begin{equation}
\mbox{minimize }S(x)+\delta(x;\Omega_0),x\in X, \label{p_1}%
\end{equation}
where $\delta(\cdot;\O_0)$ is the indicator function associated with $\O_0$ given by
\begin{equation*}\label{indicator}
\delta(x;\O_0):=\begin{cases}
0 &\text{if }\;x\in\O_0, \\
\infty & \text{otherwise}.
\end{cases}
\end{equation*}
Since each function $T_{v_{i}}(\cdot;\Omega_{i})$ is Lipschitz continuous at
$\bar{x}$, using the limiting subdifferential sum rule from \cite[Theorem
3.36]{mor}, one has
\[
0\in\partial\lbrack S(\cdot)+\delta(\cdot;\Omega_0)](\bar{x})=\sum_{i=1}%
^{n}\partial T_{v_{i}}(\bar{x};\Omega_{i})+N(\bar{x};\Omega_0)
\]
Then there exist $x_{i}^{\ast}\in\partial T_{v_{i}}(\bar{x};\Omega_{i})$ such
that
\[
-\sum_{i=1}^{n}x_{i}^{\ast}\in N(\bar{x};\Omega_0).
\]
Finally, the $x_{i}^{\ast}$ for $i=1,\ldots,n$ satisfy (1) and (2) by Theorems
\ref{limiting1} and \ref{limiting 2}. $\hfill\square$ \vspace*{0.05in}

In the convex case, we are able to obtain necessary and sufficient optimality
conditions under less restrictive assumption as in the theorem below.

\begin{Theorem} Let $X$ be a normed linear space and let $\O_i$ for $i=1,
\ldots,n$ and $\O_0$ be convex sets. Consider the optimization
problem {\rm (\ref{p})}. Suppose that there exist an element $u\in
\mbox{\rm dom }S \cap \O_0$ at which all functions
$T_{v_i}(\cdot;\O_i)$ are continuous. If $\ox$ is an optimal
solution of the problem, then there exist $x^*_i
\in X^*$ with the following properties:\\
{\rm (1)} $x^*_i\in N(\tilde{x}_i; \O_i)$, where $\tilde{x}_i:=
\Pi(\ox;\O_i)$ for $i=1,\ldots,n.$\\
{\rm (2)} $\la x_i^*, -v_i\ra \leq 1$ for all $i\in I(\ox)$ and
$\la x_i^*, v_j\ra =-1$ for all $j\in J(\ox)$.\\
{\rm (3)} $-\sum_{i=1}^n x^*_i\in N(\ox; \O_0)$. Moreover, if
$\ox\in \O_0$ satisfies {\rm(1)}, {\rm (2)}, and {\rm(3)}, then
$\ox$ is an optimal solution of the problem.
\end{Theorem}
Proof: Since $S$ is a convex function and $\Omega_0$ is a convex
set, following the proof of Theorem \ref{op1}, one sees that
$\bar{x}\in \Omega_0$ is an optimal solution of the optimization
problem (\ref{p}) if and only if
\begin{align*}
0\in &  \partial\lbrack S(\cdot)+\delta(\cdot;\Omega_0)](\bar{x})=\sum_{i=1}%
^{n}\partial T_{v_{i}}(\bar{x};\Omega_{i})+N(\bar{x};\Omega_0).
\end{align*}
Then we use the well-known convex subdifferential sum rule and Theorem
\ref{cv1} to complete the proof. $\hfill\square$

\begin{Remark} {\rm(1) With the available subdifferential formulas for
directional minimal time functions in the convex case from the paper,
we are able to develop a numerical algorithm of sub\-gradient type to
solve problem (\ref{p}) when the sets and the directions involved are
of particular shapes. See \cite{bert} for more details on the theory
of the sub\-gradient method.\\[1ex]
(2) Similar methods can be applied to the generalized Sylvester
smallest enclosing ball problem stated as follows:  given a finite
number of nonempty closed target sets $\O_i$ for $i=1,\ldots,n$ and
$n$ nonzero vectors $v_i$ for $i=1,\ldots,n$, and a nonempty closed
constraint set $\O_0$, find a point $\ox\in\O_0$ to place the
initial points of the vectors such that the vectors can reach all
the targets in the shortest time. This problem can be modeled as
follows:
\begin{equation*}\label{sv}
\mbox{\rm minimize }\max \{T_{v_i}(x;\O_i)\mid  i=1,\ldots,n\}\mbox{
subject to }x\in \O_0.
\end{equation*}}
\end{Remark}

\end{document}